\documentclass[11pt]{amsart}

\usepackage{latexsym}
\usepackage{amsmath}
\usepackage{amsfonts}
\usepackage{amssymb}
\newtheorem{theorem}{Theorem}[section]
\newtheorem{lemma}[theorem]{Lemma}
\newtheorem{proposition}[theorem]{Proposition}
\newtheorem{corollary}[theorem]{Corollary}
\newtheorem{definition}[theorem]{Definition}

\newtheorem{remark}[theorem]{Remark}

\newcommand{\prend}{ $\diamondsuit $\hfill \bigskip}

\def\Alg{\mathop{\rm Alg}}
\def\Lat{\mathop{\rm Lat}}

\def\Map{\mathop{\rm Map}}

\def\nph{\varphi}
\def\h{{\mathcal {H}}}
\def\hone{{\h_{1}}}
\def\htwo{{\h_{2}}}

\def\lba{{\mathcal {B}(\mathcal {H}_{1}, \mathcal {H}_{2})}}

\def\alga{{\mathcal {A}}}
\def\b{{\mathcal {B}}}

\def\k{{\mathcal {K}}}
\def\m{{\mathcal {M}}}

\def\t{{\mathcal {T}}}

\def\u{{\mathcal {U}}}

\def\s{{\mathcal {S}}}

\def\x{{\mathcal {X}}}
\def\y{{\mathcal {Y}}}
\def\p{{\mathcal {P}}}
\def\pone{{\mathcal {P}_{1}}}
\def\ptwo{{\mathcal {P}_{2}}}

\begin{document}

\title{Normalizers of operator algebras and reflexivity}

\author{A. Katavolos}

\address{Department of Mathematics, University of Athens,
157 84 Athens, Greece}
\email{akatavol@eudoxos.math.uoa.gr}

\author{I.G. Todorov}

\address{{\tt current address (I.G. Todorov)}
Department of Pure Mathematics, Queen's University Belfast,
Belfast BT7 1NN Northern Ireland, United Kingdom}
\email{i.todorov@qub.ac.uk}

\thanks{The second-named author was supported by a grant from
the Greek National Scholarship Foundation}



\begin{abstract}
The set of \emph{normalizers} between von Neumann (or, more generally,
reflexive) algebras $\mathcal{A}$ and $\mathcal{B}$, (that is, the set of all
operators $T$ such that $T\mathcal{A}T^{\ast}\subseteq\mathcal{B}$ and
$T^{\ast}\mathcal{B}T\subseteq\mathcal{A}$) possesses `local linear
structure': it is a union of reflexive linear spaces. These spaces belong to
the interesting class of \emph{normalizing} linear spaces, namely, those
linear spaces $\mathcal{U}$ satisfying $\mathcal{UU}^{\ast}\mathcal{U}%
\subseteq\mathcal{U}$. Such a space is reflexive whenever it is ultraweakly
closed, and then it is of the form $\mathcal{U}=\{T:TL=\phi(L)T$ for all
$L\in\mathcal{L}\}$ where $\mathcal{L}$ is a set of projections and $\phi$ a
certain map defined on $\mathcal{L}$. A normalizing space consists of
normalizers between appropriate von Neumann algebras $\mathcal{A}$ and
$\mathcal{B}$. Necessary and sufficient conditions are found for a normalizing
space to consist of normalizers between two reflexive algebras.
Normalizing
spaces which are bimodules over maximal abelian selfadjoint algebras
consist of operators `supported' on sets of the form $[f=g]$ where $f$
and $g$ are appropriate Borel functions.
They also satisfy spectral synthesis in the sense of Arveson.

\end{abstract}

\maketitle

\section{Introduction and preliminaries}

This paper is devoted to the study of a class of linear spaces of operators on
Hilbert space  which arise as normalizers of operator algebras; we call them
\emph{normalizing} subspaces.

Normalizers of von Neumann algebras (in particular, of maximal abelian
selfadjoint algebras - masas for short) are known to play an important role in
various contexts. The non-selfadjoint generalizations of  von Neumann algebras
are the\emph{ reflexive algebras} first introduced by Halmos in \cite{hal}.
We will primarily be concerned with normalizers of such algebras.
Recall that an operator $T$ normalizes a reflexive (not necessarily
selfadjoint) algebra $\mathcal{A}$ if $T$ and its adjoint satisfy
$T^{\ast}\mathcal{A}T\subseteq\mathcal{A}$.
The set of all normalizers of $\mathcal{A}$ is of course
not a linear space in general.
However, it turns out that the action of a normalizer $T$ and its adjoint
$T^{\ast}$ on the invariant projections of $\mathcal{A}$ defines a linear
space of operators all of which normalize $\mathcal{A}$. Moreover, this space is
\emph{reflexive }in the terminology of Loginov-Shulman \cite{ls} and Erdos
\cite{opmap} and is closed under the `triple product' $AB^{\ast}C$.
Thus the set of normalizers of a reflexive (not necessarily selfadjoint) algebra
(and indeed, the set of \emph{semi-normalizers} between two reflexive
algebras - see section 5) appears as the union of reflexive linear spaces which have
additional algebraic structure.
It is the interplay between linearity and normalization
that forms the subject matter of the present work.
For example,
we show that every normalizer is the norm-limit of linear combinations
of normalizing partial isometries, and every compact normalizer is the limit of
finite rank normalizers. We also show that the sum of two normalizers of
CSL algebras is again a normalizer only when both are contained in a single
reflexive masa bimodule consisting of
normalizers, and obtain generalizations of the results of
Coates \cite{coates} on normalizers of nest algebras.

These observations lead us to introduce the class of
\emph{normalizing} spaces, namely, those linear spaces of
operators which are closed under the `triple product'
$AB^{\ast}C$. These subspaces are interesting in their own right;
they generalize selfadjoint algebras of operators, and share many
properties in common with such algebras. For instance, they
satisfy an analogue of the bicommutant theorem: they are reflexive
whenever they are ultraweakly closed (unlike general
non-selfadjoint algebras). Alternatively given a set $\mathcal{L}$
of projections and a map $\phi$ defined on $\mathcal{L}$, the set
of  all operators $T$ satisfying $TL=\phi(L)T$ for all
$L\in\mathcal{L}$ is a normalizing subspace, and  every
ultraweakly closed normalizing subspace is of  this form.

A normalizing subspace $\mathcal{U}$ is a bimodule over the  selfadjoint
algebras $[\mathcal{U}^{\ast}\mathcal{U}]$ and $[\mathcal{UU}^{\ast}]$  and
induces a complete lattice isomorphism $\chi$ between the invariant
projections of the `non-degenerate parts' of these algebras. Also,
$\mathcal{U}$ normalizes the first algebra into the second (that is, $T
(\mathcal{U}^{\ast}\mathcal{U})T^{\ast}\subseteq\mathcal{UU}^{\ast}$
and $T^{\ast}(\u\u^{\ast})T\subseteq \u^{\ast}\u$ for each
$T\in\mathcal{U}$).   Conversely, we are able to characterize when
$\mathcal{U}$ normalizes a pair of reflexive (not necessarily selfadjoint)
algebras $\mathcal{A}$ and $\mathcal{B}$. Apart from the obvious relations
$\u^*\u \subseteq  \mathcal{A}$ and $\u\u^* \subseteq \mathcal{B}$,
the map $\chi$ must induce a bijection between the invariant
projection lattices of the `non-degenerate parts' of these algebras.

Thus, if $\u$ is a normalizing space then the non-degenerate parts
of the von Neumann algebras generated by $[\u^* \u]$ and $[\u \u^*]$
are Morita equivalent in the sense of Rieffel \cite{rie}.
Conversely, if $\alga$ and $\b$ are Morita equivalent W*-algebras,
then there are faithful representations of $\alga$ and $\b$ such that
the bimodule which establishes the equivalence is represented
as a normalizing space $\u$ of operators between the respective Hilbert spaces.
In this paper our concern is not with the notion of Morita equivalence
of (abstract) W*-algebras, but rather with the properties of normalizers
between (concrete) reflexive algebras and especially with the interplay
between normalizers and reflexivity.
Notice, however, that this connection between normalizers and Morita
equivalence might not have been observed had we considered normalizers
of a single algebra.

We prove that normalizing subspaces which are bimodules over two maximal
abelian selfadjoint algebras consists of operators `supported' on sets
of the form $[f=g]$ where $f$
and $g$ are appropriate Borel functions. This includes the case of normalizing
subspaces which are generated by rank one operators. In case one of the
algebras $[\mathcal{U}^{\ast}\mathcal{U}]$, $[\mathcal{UU}^{\ast}]$ is
abelian, the support of $\mathcal{U}$ turns out to be the `graph'
or the \lq reverse graph' of a Borel function. We also show that normalizing
masa-bimodules satisfy spectral synthesis in the sense of Arveson \cite{a}.
This gives a clear geometric description of the normalizers of
a CSL algebra in terms of generalized graphs
or partial graphs. These partial graphs are analogous to the ones appearing
in the work of Feldman and Moore \cite{fm} and others. In these papers,
only partial isometries normalizing
certain Cartan masas are considered, while in our work the emphasis is on the whole
reflexive linear space generated by each generalised graph. Also, we deal
with arbitrary (nonabelian and non-selfadjoint) CSL algebras.

\bigskip

The notation we use is standard; see for example \cite{dav}.
We review some definitions and facts from \cite{opmap} and \cite{ls}. Let
${{\mathcal{H}}_{1}}$ and ${{\mathcal{H}}_{2}}$ be complex Hilbert spaces,
$\mathcal{P}_{i}$ the lattice of all (orthogonal) projections on
${\mathcal{H}}_{i}$, $i=1,2$. We let $\mathcal{M(P}_{1},\mathcal{P}_{2})$
denote the set of all maps $\varphi: \mathcal{P}_{1} \to\mathcal{P}_{2}$ which
are 0-preserving and $\vee$-continuous (i.e. preserve arbitrary suprema).
Erdos \cite{opmap} shows that each
$\varphi\in\mathcal{M(P}_{1},\mathcal{P}_{2})$ uniquely defines semi-lattices
$\s_{1 \nph}\subseteq \pone$ and $\s_{2 \nph}=\nph(\pone)\subseteq \ptwo$ such that
$\nph$ is a bijection between $\s_{1 \nph}$ and $\s_{2 \nph}$ and is uniquely
determined by its restriction to $\s_{1 \nph}$.
Moreover,
${\mathcal{S}}_{1 \nph}$ is meet-complete and contains the identity projection
while ${\mathcal{S}}_{2 \nph}$ is join-complete and contains the zero projection.

Note that the set
\[
\mathop{\rm Op}\varphi=\{T\in\mathcal{B}({{\mathcal{H}}_{1}},{{\mathcal{H}%
}_{2}}):\varphi(P)^{\bot}T(P)=0\ \mbox{for each }P\in\mathcal{P}_{1}\}
\]
is also uniquely determined by $\varphi|_{\mathcal{S}_{1\varphi}}$: if $T$
satisfies $\varphi(P)^{\bot}T(P)=0$ for each $P\in\mathcal{S}_{1\varphi}$,
then $T\in\mathop{\rm Op}\varphi$.

Given a subspace $\mathcal{U\subseteq B}({{\mathcal{H}}_{1}},{{\mathcal{H}%
}_{2}})$, we define its \textbf{map} $\mathop{\rm Map} \mathcal{U}%
:\mathcal{P}_{1}\rightarrow\mathcal{P}_{2}$ by
\[
( \mathop{\rm Map} \mathcal{U})(P)= \overline{[\mathcal{U}(P)]}\quad
(P\in\mathcal{P}_{1})
\]
(where, here and in the sequel, the symbol $\overline{[\mathcal{U}(P)]}$ will
stand for the projection onto the closed subspace spanned by
$\{Sx:x\in P({{\mathcal{H}}_{1}}),S\in\mathcal{U}\}$).

If $\nph^{\ast}=\Map\u^{\ast}$, then
$\s_{1 \nph}=\{P^{\perp} : P\in\nph^* (\ptwo)\}$ \cite{opmap}.

The \textbf{reflexive hull} $\mathop{\rm Ref}\mathcal{U}$ of $\mathcal{U}$ is
defined to be the space
\[
\mathop{\rm Ref}\mathcal{U}=\{T\in{\mathcal{B}(\mathcal{H}_{1},\mathcal{H}%
_{2})}:Tx\in\overline{{\mathcal{U}}x},\ \mbox{for
each}\ x\in{{\mathcal{H}}_{1}}\}
\]
\cite{ls}. A subspace $\mathcal{U}$ is called \textbf{reflexive}
if $\mathcal{U}=\mathop{\rm Ref}\mathcal{U}$. It is easily seen
that $\mathop{\rm Ref}\mathcal{U}=\mathop{\rm Op}\mathop{\rm
Map}\mathcal{U}$ \cite{opmap}. A \emph{unital algebra} $\mathcal{A
\subseteq B(H)}$ is reflexive if and only if it is of the form
$\alga = \Alg \mathcal{L} = \{ A : L^{\bot}AL = 0 \mbox{ for all }
L \in \mathcal{L} \}$, where $\mathcal{L} = \Lat \alga$ is the
complete lattice of all invariant projections of $\alga$.

Given $\varphi\in\mathcal{M(P}_{1},\mathcal{P}_{2})$ the subspace
\[
\mathcal{V}=\{T\in\mathcal{B}({{\mathcal{H}}_{1}},{{\mathcal{H}}_{2}}):
\varphi(L) T(L^{\bot}) =0 \ \mbox{for each }L\in\mathcal{S}_{1\varphi}\}
\]
is clearly reflexive. We denote its map by $\varphi^{\bot}$. Thus
$\mathcal{V}= \mathop{\rm Op} \varphi^{\bot}$ and $\varphi^{\bot}%
\in\mathcal{M(P}_{1},\mathcal{P}_{2})$ satisfies $\varphi^{\bot}(L^{\bot}%
)\leq\varphi(L)^{\bot}$ for each $L\in\mathcal{S}_{1\varphi}$. \medskip

The following simple observations, whose proofs are routine,
will be used repeatedly.

\begin{lemma}
\label{le1}Let $\u \subseteq \lba$
be a subspace and $\mathcal{A}_{i}\subseteq\mathcal{B(H}_i)\ (i=1,2)$
be unital algebras such that
$\mathcal{A}_{2}\mathcal{UA}_{1}\subseteq\mathcal{U}$.
If $\varphi=\mathop{\rm Map}\mathcal{U}$, then
$\mathcal{S}_{i\varphi}\subseteq\mathop{\rm Lat}\mathcal{A}_{i}$. Thus
$\mathcal{A}_{1}\subseteq\mathop{\rm Alg}\mathcal{S}_{1\varphi}$ and
${\mathcal{A}}_{2}\subseteq \mathop{\rm Alg} {\mathcal{S}}_{2\varphi}$. If
additionally {$\mathcal{U}$} is reflexive, then the algebra
$\mathop{\rm Alg}{\mathcal{S}}_{1\varphi}$
(resp. $\mathop{\rm Alg} {\mathcal{S}}_{2\varphi}$) is the
largest algebra over which ${\mathcal{U}}$ is a right (resp. left) module.
\end{lemma}

\section{Normalizing spaces of operators}

The notion of  reflexivity for subspaces generalizes the corresponding  notion
defined by Halmos \cite{hal} for unital algebras. Among reflexive algebras, the
selfadjoint ones, namely the von\ Neumann algebras, have of course a
distinguished place. Note that a unital algebra $\mathcal{A}$ is selfadjoint
if and only if $\mathcal{AA}^{\ast}\mathcal{A}\subseteq\mathcal{A}$. As the
results of this paper show, the generalization of this property
to subspaces is particularly fruitful.

\begin{definition}
A subspace $\mathcal{U\subseteq B}({{\mathcal{H}}_{1}},{{\mathcal{H}}_{2}})$
is said to be \textbf{normalizing} if it is closed under the `triple product'
$(A,B,C)\rightarrow AB^{\ast}C$.
\end{definition}

\begin{remark}
\label{rem1}For a subspace $\mathcal{U\subseteq B}({{\mathcal{H}}_{1}%
},{{\mathcal{H}}_{2}})$, the following are equivalent:

\textbf{(i)} $\mathcal{U}$ is normalizing.

\textbf{(ii)} There is a unital *-algebra $\mathcal{A}_{1}\subseteq
\mathcal{B}({{\mathcal{H}}_{1}})$ such that $\mathcal{UA}_{1}\subseteq
\mathcal{U}$ and $\mathcal{U}^{\ast}\mathcal{U}\subseteq\mathcal{A}_{1}$.

\textbf{(iii)} There is a subspace $\mathcal{A}_{1}\subseteq\mathcal{B}%
({{\mathcal{H}}_{1}})$ such that $\mathcal{UA}_{1}\subseteq\mathcal{U}$ and
$\mathcal{U}^{\ast}\mathcal{U}\subseteq\mathcal{A}_{1}$.

\textbf{(iv)} There is a unital *-algebra $\mathcal{A}_{2}\subseteq
\mathcal{B}({{\mathcal{H}}_{2}})$ such that $\mathcal{A}_{2}\mathcal{U}%
\subseteq\mathcal{U}$ and $\mathcal{UU}^{\ast}\subseteq\mathcal{A}_{2}$.

\textbf{(v)} There is a subspace $\mathcal{A}_{2}\subseteq\mathcal{B}%
({{\mathcal{H}}_{2}})$ such that $\mathcal{A}_{2}\mathcal{U}\subseteq
\mathcal{U}$ and $\mathcal{UU}^{\ast}\subseteq\mathcal{A}_{2}$.
\end{remark}

\noindent\textit{Proof. } (i) $\Rightarrow$(ii) Let $\mathcal{A}_{1}%
\subseteq\mathcal{B}({{\mathcal{H}}_{1}})$ be the linear span of
$\mathcal{U}_{1}=\{S^{\ast}T:S,T\in\mathcal{U}\}\cup\{I\}$. Since
$\mathcal{UU}^{\ast}\mathcal{U\subseteq U}$, one verifies that $\mathcal{U}%
_{1}$ is a unital *-semigroup and so $\mathcal{A}_{1}$ is a unital
*-subalgebra. The properties $\mathcal{U}^{\ast}\mathcal{U\subseteq A}_{1}$
and $\mathcal{UA}_{1}\subseteq\mathcal{U}$ are immediate.

(ii)$\Rightarrow$(iii) is trivial.

(iii)$\Rightarrow$(i) If $S,T,R\in\mathcal{U}$ then $T^{\ast}R\in
\mathcal{A}_{1}$ and hence $S(T^{\ast}R)\in\mathcal{U}$.

The implications (i)$\Rightarrow$(iv)$\Rightarrow$(v)$\Rightarrow$(i) are
equally easy. $\diamondsuit$\hfill\bigskip

\noindent{\bf Remark}
Let $\u$ be a normalizing subspace, and write $\alga$ and $\b$ for the algebras
generated by $\u^* \u$ and $\u \u^*$ respectively. Then $\u$
\emph{normalizes $\b$ into
$\alga$}, in the sense that $T^* \b T \subseteq \alga$ and
$T \alga T^* \subseteq \b$ for all $T \in \u$. Conversely,
if $\u$ is a subspace of operators, which normalizes an algebra into another,
then, as we shall prove (see Proposition \ref{prop100}), $\u$ is contained
in a normalizing space of operators.

\medskip

If $\mathcal{A}$ is a unital selfadjoint algebra, then its invariant subspace
lattice $\mathcal{L}$ is orthocomplemented; thus its map
$\chi= \mathop{\rm Map}(\mathcal{A})$ (namely, the identity map of
$\mathcal{L}$) preserves orthogonality. This property characterizes maps of
normalizing subspaces.

\begin{definition}
A map $\chi\in\mathcal{M(P}_{1},\mathcal{P}_{2})$ is said to be an
\textbf{ortho-map  }if $\chi(L)\bot\chi(L^{\bot})$ for each $L\in
\mathcal{S}_{1\chi}$.
\end{definition}

The following theorem shows the connection between ortho-maps and
normalizing spaces. Statement (b)(ii) corresponds to the von Neumann Bicommutant
Theorem.

\begin{theorem}
\label{th0.1}\label{th0} \label{24}
\textbf{(a)} Let $\varphi\in\mathcal{M(P}_{1},\mathcal{P}_{2})$ and
\[
\mathcal{U}=\{T\in\mathcal{B}({{\mathcal{H}}_{1}},{{\mathcal{H}}_{2}%
}):TL=\varphi(L)T\ \mbox
{for each }L\in\mathcal{S}_{1\varphi}\}.
\]
Then $\mathcal{U\ =}\mathop{\rm Op}\varphi\cap\mathop{\rm Op}\varphi^{\bot}$
is a normalizing space.

\textbf{(b)} Let $\mathcal{U\subseteq B}({{\mathcal{H}}_{1}},{{\mathcal{H}%
}_{2}})$ be a normalizing subspace. Then

(i) $\mathop{\rm Map}\mathcal{U}$ is an ortho-map and

(ii) $\mathop{\rm Ref}(\mathcal{U)}=\mbox{cl}_{WOT}(\mathcal{U)}=\mbox
{cl}_{uw}(\mathcal{U)}$.

\textbf{(c)} Let $\chi\in\mathcal{M(P}_{1},\mathcal{P}_{2})$ be an ortho-map.
Then
\[
\mathop{\rm Op}\chi=\{T\in\mathcal{B}({{\mathcal{H}}_{1}},{{\mathcal{H}}_{2}%
}):TL=\chi(L)T\ \mbox{for
each }L\in\mathcal{S}_{1\chi}\}.
\]
\end{theorem}

\noindent\textit{Proof. } (a) The equality
$\u = \mathop{\rm Op} \varphi\cap\mathop{\rm Op} \varphi^{\bot}$
is easily verified.

Let $S,T,R\in \u$. Then for each $L\in\mathcal{S}_{1\varphi}$ we have
\[
(ST^{\ast}R)L=ST^{\ast}\varphi(L)R=SLT^{\ast}R=\varphi(L)(ST^{\ast}R)
\]
and so $ST^{\ast}R\in\u$ which shows that $\u$ is normalizing.

\bigskip(b) Let $\mathcal{A}_{i}\subseteq\mathcal{B}({\mathcal{H}}_{i})$ be
unital *-algebras such that $\mathcal{U}^{\ast}\mathcal{U\subseteq A}_{1}$,
$\mathcal{UU}^{\ast}\mathcal{\subseteq A}_{2}$ and $\mathcal{A}_{2}%
\mathcal{UA}_{1}\subseteq\mathcal{U}$ (Remark \ref{rem1}).

\smallskip(i) Let $\chi= \mathop{\rm Map} \mathcal{U}$.
We show that $\chi(L)\bot\chi(L^{\bot})$
for each $L\in\mathop{\rm Lat} \mathcal{A}_{1}$ (this will suffice since
$\mathcal{S}_{1\chi}\subseteq\mathop{\rm Lat}\mathcal{A}_{1}$
by Lemma \ref{le1}). Indeed,
for each $S,T\in\mathcal{U}$ we have $S^{\ast}T\in \mathcal{A}_{1}$ and so,
if $\xi\in L$ and $\eta\in L^{\bot}$ then
\[
\langle T\xi,S\eta\rangle=\langle S^{\ast}T\xi,\eta\rangle=0.
\]
This shows that $T\xi\bot\chi(L^{\bot})$ for each
$T\in\mathcal{U}$ and hence $\chi(L)\bot\chi(L^{\bot})$.

\bigskip(ii)
The properties $\u^{\ast}\u\subseteq\alga_1, \u\u^{\ast}\subseteq\alga_2$
and $\alga_2\u\alga_1\subseteq\u$ ensure that the set
\[
\mathcal{C}=\left(
\begin{array}
[c]{cc}%
\mathcal{A}_{2} & \mathcal{U}\\
\mathcal{U}^{\ast} & \mathcal{A}_{1}%
\end{array}
\right)  =\left\{  \left(
\begin{array}
[c]{cc}%
B & T\\
S^{\ast} & A
\end{array}
\right)  :A\in\mathcal{A}_{1},B\in\mathcal{A}_{2},S,T\in\mathcal{U}\right\}
\]
is a unital *-subalgebra of $\mathcal{B}({{\mathcal{H}}_{2}}\oplus
{{\mathcal{H}}_{1}})$ and so $\mbox{cl}_{WOT}(\mathcal{C})=\mbox{cl}%
_{uw}(\mathcal{C}) = {\mathcal{C}}^{\prime\prime}$ by the von Neumann
bicommutant theorem. But, since ${\mathcal{C}}$ is a unital *-algebra, it is
easy to verify that ${\mathcal{C}}^{\prime\prime} = \mathop{\rm Alg}
\mathop{\rm Lat}{\mathcal{C}} = \mathop{\rm Ref}(\mathcal{C)}$.
This implies in particular that $\mbox{cl}_{WOT}(\mathcal{U})=
\mbox{cl}_{uw}(\mathcal{U}) =\mathop{\rm Ref}(\mathcal{U)}$.

\bigskip(c) If $T\in\mathcal{B}({{\mathcal{H}}_{1}},{{\mathcal{H}}_{2}})$
satisfies $TL=\chi(L)T$ for each $L\in\mathcal{S}_{1\chi}$, then
$TL=\chi(L)TL$ and so $T\in\mathop{\rm Op}\chi$ since the latter is determined
by $\mathcal{S}_{1\chi}$. If $\chi$ is an ortho-map and
$T\in\mathop{\rm Op}\chi$, then for each $L\in\mathcal{S}_{1\chi}$ the
relation $TL^{\bot}%
=\chi(L^{\bot})TL^{\bot}$ gives $\chi(L)TL^{\bot}=\chi(L)\chi(L^{\bot
})TL^{\bot}=0$ since $\chi(L)\bot\chi(L^{\bot})$. Adding to this the relation
$TL=\chi(L)TL$ gives $TL=\chi(L)T$ as required. $\diamondsuit$\hfill\bigskip

We isolate two consequences of this theorem for emphasis.

\begin{corollary}
\textbf{(i)} The w*-closure of a normalizing subspace is reflexive and coincides
with its WOT-closure.

\textbf{(ii)} A reflexive subspace is normalizing if and only if its map is an
ortho-map.
\end{corollary}

\noindent\textbf{Remarks} We do not know whether a WOT-closed subspace whose
map is an ortho-map must be normalizing. Note that the map of a unital
algebra is an
ortho-map if and only if it its invariant subspace lattice is
orthocomplemented. Thus the question, within the class of unital algebras,
reduces to the well-known reductive algebra problem \cite{rr}: must a
WOT-closed algebra whose invariant lattice is orthocomplemented be selfadjoint?

A subspace whose map is an ortho-map need not be normalizing. Indeed there exist
nonselfadjoint transitive algebras (even triangular ones - see \cite{kadsin}).
However, as we show below (Corollary \ref{c16.1}) an ultraweakly closed
subspace whose map is an ortho-map must be normalizing, provided it is a masa bimodule.

\bigskip

A crucial property of von Neumann algebras is that they are generated by
their projections. Of course, normalizing spaces need not contain any
(nontrivial) projections; their role is played by the partial isometries.

\begin{proposition}\label{prop-1}
If $\mathcal{U}$ is an ultraweakly closed normalizing space and $A=U|A|$ is the
polar decomposition of an element of $\mathcal{U}$, then $U$ is a partial
isometry in $\mathcal{U}$ and $Uf(|A|)\in\u$, for every Borel function
$f$ on the spectrum $sp(|A|)$ of $|A|$.
Moreover, $\mathcal{U}$ is the norm-closed linear span of the partial
isometries it contains.
\end{proposition}

\noindent\textit{Proof. } Let ${\mathcal{A}}_{1}$ be a von Neumann algebra
with the property ${\mathcal{U}}^{\ast}{\mathcal{U}}\subseteq {\mathcal{A}}_{1}$
and ${\mathcal{U}}{\mathcal{A}}_{1}\subseteq {\mathcal{U}}$. Then $A^{\ast}%
A\in\mathcal{A}_{1}$ and so $|A|\in\mathcal{A}_{1}$. We have $U=$w$^{\ast}$%
-$\lim_{\varepsilon\rightarrow0}A(|A|+\varepsilon)^{-1}$ and so $U\in
\mathcal{U}$.
Since $|A|\in\alga_1$, for every Borel function $f$ on $sp(|A|)$, the
operator
$f(|A|)$ is in $\alga_1$ as well and since $\u$ is a right $\alga_1$-module,
it follows that $Uf(|A|)\in\u$.

Given $A \in \u$ and $\epsilon >0$, there are spectral projections
$P_1, \ldots, P_n$ of $|A|$ and scalars $c_1, \ldots, c_n$ such that
$\| |A| - \sum c_iP_i \| < \epsilon$. Thus
$\| A - \sum c_iUP_i \| < \epsilon$ and each $UP_i$ is a partial isometry
in $\u$, since the initial projection of $U$ is the range projection of
$|A|. \qquad \diamondsuit$

\medskip

A w*-closed normalizing subspace need not contain (nonzero) finite rank
operators. We show that then it cannot contain compact operators, contrary to
the situation in general reflexive subspaces \cite{hmt}.

\begin{corollary}
\label{prop88} Let ${\mathcal{U}}$ be an ultraweakly closed normalizing space
and suppose that $K$ is a compact operator in ${\mathcal{U}}$. Then $K$ can be
approximated in the norm topology by finite rank operators in ${\mathcal{U}}$.
Moreover, if $K$ belongs to some Schatten class ${\mathcal{C}}_{p}$ then
it can be approximated by finite
rank operators in ${\mathcal{U}}$ in the $p$-norm topology.
\end{corollary}

\noindent\textit{Proof. }
Immediate from Proposition \ref{prop-1}.
\prend

Let us note that the last corollary can be inferred in a different way.
Namely, consider the ``matrix'' algebra ${\mathcal{C}}$ defined in the proof
of Theorem \ref{th0.1}. Then $\mathcal{C}$ is a von Neumann algebra (Remark
\ref{rem1}) and the result is clear for $\mathcal{C}$.

\bigskip

If $\mathcal{U}$ is a normalizing space, then its rank one subspace (i.e. the
linear span of the rank one operators contained in $\mathcal{U}$) is also
normalizing. The next proposition characterizes this subspace in terms of the
map of $\mathcal{U}$. We denote the rank one operator sending $\xi$ to
$\langle\xi,x\rangle y$ by the symbol $y\otimes x^{\ast}$.

\begin{proposition}
\label{l10} Let $\chi$ be an ortho-map. A rank one operator $y\otimes
x^{\ast}$ belongs to $\mathop{\rm Op}\chi$ if and only if, for each
$\ L\in{\mathcal{S}}_{1\chi}$,
\[
Lx\neq0\Leftrightarrow L^{\perp}x=0\Leftrightarrow\chi(L)^{\perp
}y=0\Leftrightarrow\chi(L)y\neq0.
\]
\end{proposition}

\noindent\textit{Proof. }By Theorem \ref{th0.1}, $\mathop{\rm Op}\chi
=\{T\in{\mathcal{B}(\mathcal{H}_{1},\mathcal{H}_{2})}:TL=\chi(L)T,\ L\in
{\mathcal{S}}_{1\chi}\}$. Suppose that the rank one operator $y\otimes
x^{\ast}$ belongs to $\mathop{\rm Op}\chi$. Then, for each $L\in{\mathcal{S}%
}_{1\chi}$, $(\chi(L)^{\perp}y)\otimes(Lx)^{\ast}=\chi(L)^{\perp}(y\otimes
x^{\ast})L=0$ and $(\chi(L)y)\otimes(L^{\perp}x)^{\ast}=\chi(L)(y\otimes
x^{\ast})L^{\perp}=0$. So for each $L\in{\mathcal{S}}_{1\chi}$, (either
$\chi(L)^{\perp}y=0$ or $Lx=0$) and (either $\chi(L)y=0$ or $L^{\perp}x=0$).
If $Lx\neq0$ and $L^{\perp}x\neq0$, then by the above conditions, we conclude
that $\chi(L)^{\perp}y=\chi(L)y=0$, so $y=0$, which is impossible. Because
$Lx$ and $L^{\perp}x$ cannot both be zero, we conclude that $Lx\neq
0\Leftrightarrow L^{\perp}x=0$. In the same way $\chi(L)^{\perp}%
y\neq0\Leftrightarrow\chi(L)y=0$. Similarly, if $Lx\neq0$, then
$\chi(L)^{\perp}y=0$ and, conversely, if $\chi(L)^{\perp}y=0$, then
$\chi(L)y\neq0$ and so $L^{\perp}x=0$. The converse is trivial.
$\diamondsuit $\hfill\bigskip

We would now like to show how the general theory of Erdos \cite{opmap}
specializes in the case of ortho-maps.

\begin{definition}\label{dess}
For a map $\chi$, set $I_-=\chi(I)\in\ptwo$ and $0_+=\vee\{P :
\chi(P)=0\}\in \pone$. The map $\chi$ will be called
\textbf{essential} if $0_+=0$ and $I_-=I$. A subspace of operators
$\u$ will be called \textbf{essential}, if $\Map\u$ is essential.
\end{definition}

\begin{theorem}
\label{th0.001} \label{210}
Suppose that $\u$ is a normalizing space of operators and
let $\chi=\Map\u$. Then the semi-lattices
${\mathcal{S}}_{1}$ and ${\mathcal{S}}_{2}$ of $\chi$ are complete
ortho-lattices and $\chi$ is a complete ortho-lattice isomorphism of $\s_1$
onto $\s_2$. Moreover, if
$\t_1 =\s_1|_{0_+^{\perp}\hone}$ and $\t_2=\s_2|_{I_{-}\htwo}$, then
$\t_1$ and $\t_2$ are the projection lattices of von Neumann algebras and

(i) $\mathop{\rm Alg}\mathcal{S}_{1}=\left\{\left(
\begin{array}
[c]{cc}%
A & 0\\
C & D
\end{array}
\right)  :\ A\in{\mathcal{T}_{1}}^{\prime},C\in{\mathcal{B}}(0_{+}^{\perp
}{{\mathcal{H}}_{1}},0_{+}{{\mathcal{H}}_{1}}),D\in{\mathcal{B}}%
(0_{+}{{\mathcal{H}}_{1}},0_{+}{{\mathcal{H}}_{1}})\right\}$ and
$({\mathcal{U}}^{\ast}{\mathcal{U}})^{\prime\prime}=
\left\{\left(
\begin{array}
[c]{cc}%
A & 0\\
0 & \lambda I
\end{array}
\right)  :A\in{\mathcal{T}_{1}}^{\prime}, \ \lambda \in \mathbb{C} \right\}$

(ii) $\mathop{\rm Alg}\mathcal{S}_{2}=\left\{\left(
\begin{array}
[c]{cc}%
B & C\\
0 & D
\end{array}
\right)  :B\in{\mathcal{T}_{2}}^{\prime},C\in{\mathcal{B}}(I_{-}^{\perp
}{{\mathcal{H}}_{2}},I_{-}{{\mathcal{H}}_{2}}),D\in{\mathcal{B}}(I_{-}^{\perp
}{{\mathcal{H}}_{2}},I_{-}^{\perp}{{\mathcal{H}}_{2}})\right\}$ and
$({\mathcal{U}}{\mathcal{U}}^{\ast})^{\prime\prime}=\left\{\left(
\begin{array}
[c]{cc}%
B & 0\\
0 & \lambda I
\end{array}
\right)  :B\in{\mathcal{T}_{2}}^{\prime}, \ \lambda \in \mathbb{C} \right\}$.

(iii) $\Lat(\u^{\ast}\u)=\{L\in\pone : \chi(L)\perp\chi(L^{\perp})\} =
\{L_1 \oplus L_2 : L_1 \in \mathcal{T}_{1} \}$.
\end{theorem}
\proof
Put $\chi^{\ast}=\Map\u^{\ast}$. Recall that
$\s_1=\{P^{\perp} : P\in\chi^{\ast}(\ptwo)\}$ and $\s_2=\chi(\pone)$.
First assume that $\u$ is essential. Let $\x=\Lat(\u^{\ast}\u)$ and
$\y=\Lat(\u\u^{\ast})$. Since $\u^{\ast}\u$ and $\u\u^{\ast}$ are
selfadjoint sets of operators, $\x$ and $\y$ are complete orthocomplemented
projection lattices.

We will show that $\chi$ is an ortho-isomorphism from
$\x$ onto $\y$.

First, if $L\in\x$, then given $T,S\in\u$ we have
$(TS^{\ast})(R(L\hone))=(TS^{\ast}R)(L\hone)\subseteq\chi(L)\htwo$ for all
$R\in\mathcal{U}$.
Hence $(TS^{\ast})(\chi(L)\htwo)\subseteq\chi(L)\htwo$ and so
$\chi(L)\in\mathcal{Y}$.
Similarly, $\chi^{\ast}(\y) \subseteq \x$.

Next, we show that, if $L\in\x$, then $\chi^{\ast}\chi(L)=L$.
For $T,S\in\mathcal{U}$ we have $(T^{\ast}S)L=L(T^{\ast}S)$ since
$L\in(\u^*\u)'$ and
so $L^{\perp}(T^{\ast}S)L=0$ and $L(T^{\ast}S)L^{\bot}=0$. The
first relation gives $L^{\perp}T^{\ast}(\chi(L))=0$ for each $T\in\u$
and so $\chi^{\ast}\chi(L)\leq L$. Similarly the second gives
$\chi^{\ast}\chi(L^{\bot})\leq L^{\perp}$.
But $\chi^{\ast}\chi(L\vee L^{\bot})=\chi^{\ast}\chi(I)=I$ and so
$\chi^{\ast}\chi(L)\vee\chi^{\ast}\chi(L^{\bot})=I$.
The relation $\chi^*\chi(L)=L$ now follows readily.
Applying the same arguments to the space $\u^{\ast}$, it follows that,
if $M\in\y$, then $\chi\chi^{\ast}(M)=M$.

So we have shown that $\chi$ maps $\mathcal{X}$ bijectively to $\mathcal{Y}$
and its inverse is $\chi^{\ast}$. But $\chi$ preserves orthogonality (indeed
if $T,S\in\mathcal{U}$ and $L\in\mathcal{X}$ then $\langle TLx,SL^{\bot
}y\rangle=$ $\langle S^{\ast}TLx,L^{\bot}y\rangle=0$ since $L^{\perp}S^{\ast
}TL=0$) and since it is $\vee$-continuous and preserves $0$ and $I$, it
follows that it is $\wedge$-continuous as well and hence an ortho-isomorphism.

Finally, notice that $\s_2=\mathcal{Y}$ and $\s_1=\mathcal{X}$.
Indeed, given $P\in\pone$ and $T,S\in\mathcal{U}$ we have
$\chi(P)^{\perp}(TS^{\ast})RP=0$
for all $R\in\mathcal{U}$ (since $TS^*R\in\u$ and $\overline{[\u(P)]}= \chi(P)$),
hence $\chi(P)^{\perp}(TS^{\ast})\chi(P)=0$.
Thus $\chi(P)\in\y$, for each $P\in\pone$ and so
$\mathcal{S}_{2\chi}=\chi(\mathcal{P}
_{1})\subseteq\mathcal{Y}$. Also if $M\in\mathcal{Y}$ then $M=\chi
(L)\in\mathcal{S}_{2\chi}$ where $L=\chi^{\ast}(M)$. This shows that
$\y=\s_2$. The proof that $\x=\s_1$ is similar.

\smallskip

Now relax the assumption that $\u$ is essential.
Note that $0_{+}^{\bot}$ acts as the identity on
$\mathcal{U}^{\ast}\mathcal{U}$ and hence
$0_{+}\in(\mathcal{U}^{\ast}\mathcal{U})^{\prime\prime}$.
Similarly $I_{-}\in(\mathcal{UU}^{\ast})^{\prime\prime}$ acts as the identity
on $\mathcal{UU}^{\ast}$ and
$\u(\hone)=\u(0_{+}^{\bot}\hone)\subseteq I_{-}\htwo$. Let
$\mathcal{K}_{1}=0_{+}^{\bot}\hone$, $\mathcal{K}_{2}=I_{-}\htwo$ and
$\u_{o}=\u|_{\mathcal{K}_{1}}\subseteq\mathcal{B(K}_{1},\mathcal{K}_{2})$.
Then $\mathcal{U}_{o}$ is normalizing and essential. If $\chi_{o}$ is
the map of $\mathcal{U}_{o}$ then clearly its semilattices are
$\s_{1\chi_{o}} = \s_{1\chi}|_{\k_1} = \t_1$ and
$\s_{2\chi_{o}} = \s_{2\chi}|_{\k_2} = \t_2$. By what was shown above,
$\t_1 = \Lat (\u_{o}^{*} \u_{o})$, $\t_2 = \Lat (\u_{o} \u_{o}^{*})$ and
$\chi_o : \t_1 \to \t_2$ is an isomorphism. Also
$\mathcal{(U}^{\ast}\mathcal{U})^{\prime\prime}=
\overline{[\u_{o}^*\u_{o}]}\oplus\mathbb{C}I_{\k_1}$ and
$\mathcal{(UU}^{\ast})^{\prime\prime} =
\overline{[\u_{o}\u_{o}^*]}\oplus\mathbb{C}I_{\k_2}$.

It follows that $\Lat (\u^* \u) = \{ L_1 \oplus L_2 : L_1 \in \t_1 \}$,
while the other equality of \textit{(iii)} is readily verified, in fact
for every subspace $\u$ with map $\chi$.

Now a projection $M\in\mathcal{P}_{2}$ is in $\s_{2\chi}=\chi(\p_1)$ if and
only if it is of the form $M=M_{1}\oplus0$ with respect to the decomposition
$\mathcal{H}_{2}=\mathcal{K}_{2}\oplus\mathcal{K}_{2}^{\bot}$, where
$M_{1}\in \mathcal{S}_{2\chi_{o}}$.
Similarly working with $\mathcal{U}^{\ast}$ we obtain that
$\chi^*(\mathcal{P}_{2})=\{L_{1}\oplus0:L_{1}\in\mathcal{S}_{1\chi_{o}}\}$
with respect to the decomposition $\mathcal{H}_{1}=\mathcal{K}_{1}\oplus
\mathcal{K}_{1}^{\bot}$.  Thus
\[
\mathcal{S}_{1\chi}=\{L_{1}\oplus I:L_{1}\in\t_1\},\qquad
\mathcal{S}_{2\chi}=\{M_{1}\oplus0:M_{1}\in\t_2\}.
\]
From this we see that $X\in \Alg \mathcal{S}_{1\chi}$ if and only if
$(L_{1}\oplus I)^{\bot}X(L_{1}\oplus I)=0$ for each $L_{1}\in\t_1$,
and claim \textit{(i)} follows easily. The proof of \textit{(ii)} is
similar.
\prend

\noindent \textbf{Remarks (i) } It is not difficult to show that
the von Neumann algebras $\t_1''$ and $\t_2''$ are in fact
*-isomorphic, and so the strongly closed algebras generated by
$\u^* \u$ and $\u \u^*$ are Morita equivalent \cite{rie}. We have
preferred the direct approach above, which is sufficient for our
needs.

\textbf{(ii) }
The first part of the above theorem remains valid even if $\chi$ is
not given as a map of some normalizing subspace, but is an arbitrary
ortho-map. That is, the following lattice-theoretic result holds:
If $\chi$ is an ortho-map, then its semi-lattices $\s_1$ and $\s_2$
are in fact complete ortho-lattices and $\chi:\s_1\longrightarrow \s_2$
is a complete ortho-isomorphism \cite{it}.

\medskip

The next corollary shows that the semilattices of the map of a normalizing space
become \emph{reflexive} lattices, if
the extreme elements $0$ and $I$ are adjoined.

\begin{corollary}
\label{c0.0002} Suppose that ${\mathcal{S}}_{1}$ and ${\mathcal{S}}_{2}$ are
the (semi-)lattices of the map $\chi$ of some normalizing subspace. Then
$\Lat\Alg\s_1={\mathcal{S}}_{1}\cup\{0\}$ and
$\Lat\Alg\s_2={\mathcal{S}}_{2}\cup\{I\}$. Moreover,
in the notation of Theorem \ref{th0.001},
the orthocomplemented projection lattice generated by $\s_1$, is
${\mathcal{S}}_{1}\cup\{P_{1}\oplus0 : P_{1} \in{\mathcal{T}_{1}}\}$, while the
orthocomplemented projection lattice generated by $\s_2$ is
${\mathcal{S}}_{2} \cup\{Q_{1}\oplus I : Q_{1}\in{\mathcal{T}_{2}}\}.$
\end{corollary}

\noindent\textit{{Proof. }}
Let $\Sigma_1=\Lat\Alg\s_1$, $\Sigma_2=\Lat\Alg\s_2$ and
$\widetilde{\Sigma_1}$ and $\widetilde{\Sigma_2}$ be the orthocomplemented
lattices generated by $\Sigma_1$ and $\Sigma_2$ respectively.
Suppose that a projection $Q$ is left invariant by
the algebra $\mathop{\rm Alg}{\mathcal{S}}_{1}$. Then by Theorem
\ref{th0.001} it must be contained in $\Lat ( \u^* \u)$; thus $Q$ is of the form
$Q=Q_{1}\oplus Q_{2}$, where $Q_{1}\in \t_1$. Since
$0_{+}\mathop{\rm Alg}\s_1|_{0_{+}\hone}=\b(0_{+}\hone,0_{+}\hone)$,
it follows that either
$Q_{2}=I$ or $Q_{2}=0$. If $Q_{2}=I$ then $Q\in\mathcal{S}_{1}$. If
$Q_{2}=0$, then $CQ_{1}=0$ for each
$C\in \b(0_{+}^{\perp}\hone,0_{+}\htwo)$
and it follows that $Q_{1}=0$.
Thus we have that $\Sigma_{1}={\mathcal{S}}_{1}\cup\{0\}$.

It is easy to check that the set
$\s_1 \cup \{ Q_1 \oplus 0 : Q_1 \in \t_1 \}$ is an orthocomplemented lattice.
Since it contains $\s_1$, it must equal $\widetilde{\Sigma_{1}}$.

The other identities are proved by the same arguments. $\diamondsuit$

\bigskip

Recall that a subspace ${\mathcal{U}}$ is called \textbf{strongly reflexive}
\cite{eks}, if there is a set of rank one operators ${\mathcal{R}}$ such that
${\mathcal{U}}=\mathop{\rm Ref}{\mathcal{R}}$. We wish to describe the
strongly reflexive normalizing subspaces. The unital algebra case might be known;
we include a proof for completeness.

\begin{lemma}\label{lemma10}
Every strongly reflexive selfadjoint unital algebra contains a totally atomic
maximal abelian selfadjoint algebra (a masa).
\end{lemma}

\proof
If $\mathcal{A}\subseteq\mathcal{B(H)}$ is a strongly reflexive
selfadjoint unital algebra, then it is a von Neumann algebra, since it is
reflexive. If $\mathcal{L}=$ $\mathop{\rm Lat}{\mathcal{A}}$, then
$\mathcal{L}=\mathop{\rm Lat}{\mathcal{R}}$, where $\mathcal{R\subseteq A}$ is
the rank one subalgebra of $\mathcal{A}$,  by strong reflexivity. It follows
from \cite{lon} that $\mathcal{L}$ is completely distributive. But
$\mathcal{L}$ is orthocomplemented, hence, by Tarski's theorem (see \cite{birk},
p.119), it must be a complete atomic Boolean lattice. Thus, as is well-known,
 $\mathcal{L}$  is commutative, hence
$\mathcal{L\subseteq A}$. Therefore there exists a totally atomic masa
$\m$ with $\mathcal{L \subseteq M \subseteq A}$.
\prend

\begin{proposition}
Every strongly reflexive normalizing subspace
$\mathcal{U\subseteq B(H}_{1},\mathcal{H}_{2})$
is a bimodule over two totally atomic masas.
\end{proposition}

\proof
Let $x\in \overline{\u^*\htwo} \equiv \k_1$. If $S,T\in\mathcal{U}$, then, since
$\mathcal{U}$ is strongly reflexive,  $S^{\ast}Tx$ can be approximated by
vectors of the form $S^{\ast}Rx$, where $R$ belongs to the rank one subspace
of $\mathcal{U}$.
The strong closure of the linear span of $\u^*\u|_{\k_1}$ is
therefore a strongly reflexive von Neumann algebra acting on the space $\k_1$.
Thus it must contain a totally atomic masa by Lemma \ref{lemma10}.
Since ${\mathcal{B}}(\k_1^{\bot}, \k_1^{\bot})$
also contains a totally atomic masa, it follows from Theorem \ref{th0.001}
that $\Alg\s_1$ contains a totally atomic masa.

The proof that $\Alg\s_2$ contains a totally atomic masa is
identical. Since $\u$ is an $(\Alg\s_2,\Alg\s_1)$-bimodule, we are
done. \prend

We will in fact prove (Theorem \ref{th11})
that normalizing strongly reflexive subspaces are of the form
$\oplus_{\lambda}{\mathcal{B}}({\mathcal{H}}_{\lambda},{\mathcal{K}}_{\lambda})$,
for some mutually orthogonal families of subspaces
$\mathcal{H}_{\lambda}\subseteq\mathcal{H}_1$ and
${\mathcal{K}}_{\lambda}\subseteq
\mathcal{H}_2$.

\section{Normalizing masa bimodules}

In this and the next section, we specialize to subspaces $\mathcal{U\subseteq
B(}{{\mathcal{H}}_{1}},{{\mathcal{H}}_{2}})$ which are bimodules over maximal
abelian selfadjoint algebras (for short masas) ${\mathcal{D}_{1}\subseteq
}\mathcal{B}({{\mathcal{H}}_{1}})$ and ${\mathcal{D}_{2}\subseteq}%
\mathcal{B}({{\mathcal{H}}_{2}})$ in the sense that ${\mathcal{D}_{2}}%
${$\mathcal{U}$}${\mathcal{D}_{1}}\subseteq${$\mathcal{U}$}. As we will be
using measure-theoretic arguments, we will make the blanket assumption that
from now on \emph{ all Hilbert spaces will be separable}. If ${\mathcal{L}%
_{1}}$ and ${\mathcal{L}_{2}}$ are complete projection lattices, we let
\[
{\mathcal{M}}({\mathcal{L}_{1}},{\mathcal{L}_{2}})=\{\varphi\in{\mathcal{M}%
}({\mathcal{P}_{1}},{\mathcal{P}_{2}}):{\mathcal{S}}_{i\varphi}\subseteq
{\mathcal{L}}_{i},\ i=1,2\}.
\]

\begin{theorem}
\label{th6} A w*-closed masa-bimodule ${\mathcal{U}}$ is normalizing if and
only if there exist nests $\mathcal{N}_{1},\mathcal{N}_{2}$ and a map
$\varphi\in{\mathcal{M}}({\mathcal{N}_{1}},{\mathcal{N}_{2}})$ such that
${\mathcal{U}}=\mathop{\rm Op}(\varphi)\cap\mathop{\rm Op}(\varphi^{\bot})$.
\end{theorem}

\noindent\textit{Proof. }Every space of the form ${\mathcal{U}}=\mathop{\rm
Op}(\varphi)\cap\mathop{\rm Op}(\varphi^{\bot})$ is normalizing (Theorem
\ref{24}); and if $\varphi$ is a nest map, then $\mathcal{U}$ is a masa bimodule.

For the converse, write ${\mathcal{U}}$ as $\mathop{\rm Op}\chi$ for an
ortho-map $\chi$. Recall (Theorem \ref{210}) that the lattice $\mathcal{S}%
_{1o}=\mathcal{S}_{1\chi}|_{(0_{+}\mathcal{H}_{1})^{\bot}}$ is
orthocomplemented, hence, since $\mathcal{S}_{1\chi}\subseteq\mathcal{D}_{1}$,
a (complete) Boolean lattice. There exists a complete nest
$\mathcal{N}_o \subseteq \mathcal{S}_{1o}$
on the space $(0_{+}{{\mathcal{H}}_{1}})^{\perp}$, generating
$\mathcal{S}_{1o}$ as a complete Boolean lattice \cite{a}. Then the nest
$\mathcal{N}_1=\{ Q \oplus I :Q\in\mathcal{N}_o \}\cup\{0\}$ on
$\mathcal{H}_1$ will generate the complete Boolean lattice
$\widetilde{\Sigma_{1}}$ generated by ${\mathcal{S}}_{1\chi}$. We define the
map $\varphi$ to be the restriction of $\chi$ to ${\mathcal{N}_{1}}$. From the
fact that ${\mathcal{N}_{1}}$ is contained in $\mathcal{S}_{1\chi}\cup\{0\}$ it
is clear that the left semi-lattice of $\varphi$, ${\mathcal{S}}_{1\varphi}$,
equals ${\mathcal{N}_{1}}\setminus\{0\}$ (or ${\mathcal{N}_{1}}$, if $0_{+}%
=0$). If $X\in\mathop{\rm Op}(\chi)$, then certainly $X\in\mathop{\rm
Op}(\varphi)\cap\mathop{\rm Op}(\varphi^{\perp})$. Conversely, let
$X\in\mathop{\rm Op}(\varphi)\cap\mathop{\rm Op}(\varphi^{\perp})$. For each
$N\in{\mathcal{N}_{1}}$, we have  $XN{{\mathcal{H}}_{1}}\subseteq
\varphi(N){{\mathcal{H}}_{2}}=\chi(N){{\mathcal{H}}_{2}}$ and $XN^{\perp
}{{\mathcal{H}}_{1}}\subseteq\varphi(I)\varphi(N)^{\perp}{{\mathcal{H}}_{2}%
}=\chi(I)\chi(N)^{\perp}{{\mathcal{H}}_{2}}=\chi(N^{\perp}){{\mathcal{H}}_{2}%
}$. Because ${\mathcal{N}_{1}}$ generates $\widetilde{\Sigma_{1}}$ and
$\chi$ preserves arbitrary unions and intersections, it follows that
$X\in\mathop{\rm Op}\chi$. $\diamondsuit$\hfill\bigskip

To state the next results, we need to recall some terminology from \cite{eks}.
Let $(X,\mu)$, $(Y,\nu)$ be standard Borel spaces and let $\mathcal{D}%
_{1},\mathcal{D}_{2}$ be the multiplication masas on the corresponding $L^{2}$
spaces ${{\mathcal{H}}_{1}},{{\mathcal{H}}_{2}}$, with projections denoted by
$E(\alpha)\in\mathcal{D}_{1}$ and $F(\beta)\in\mathcal{D}_{2}$, where
$\alpha\subseteq X$ and $\beta\subseteq Y$ are Borel sets. A bounded operator
$T:{{\mathcal{H}}_{1}}\rightarrow{{\mathcal{H}}_{2}}$ is said to be
\textbf{supported} by a set $\kappa\subseteq X\times Y$ if $F(\beta
)TE(\alpha)=0$ whenever $(\alpha\times\beta)\cap\kappa=\emptyset$. It is shown in
\cite{eks} that a masa bimodule $\mathcal{U}$ is reflexive precisely when
there exists a set $\kappa\subseteq X\times Y$ such that $\mathcal{U}$
consists of all operators supported by $\kappa$. This set is uniquely defined
up to marginal equivalence, and is called the $\omega$\textbf{-support} of
$\mathcal{U}$; its complement is (marginally equivalent to)
a countable union of Borel rectangles.

\begin{theorem}
\label{p1} Let $\mathcal{H}_{1}=L^{2}(X,\mu)$, $\mathcal{H}_{2}=L^{2}(Y,\nu)$
and let $\mathcal{U\subseteq B(H}_{1},\mathcal{H}_{2})$ be a w*-closed
normalizing space which is a bimodule over the multiplication masas
$\mathcal{D}_{i}\subseteq\mathcal{B(H}_{i})$. Then there exist Borel functions
$f:X\rightarrow\lbrack0,1]$ and $g:Y\rightarrow\lbrack0,1]$ such that the
$\omega$-support of $\mathcal{U}$ is the set
$\kappa=\{(x,y)\in X\times Y:f(x)=g(y)\}$.

Conversely, the set of all operators supported by a set of the above form is a
w*-closed normalizing masa bimodule.
\end{theorem}

\noindent\textit{Proof }\ Let $\chi$ be the map of $\mathcal{U}$, and let
$X_{o}\subseteq X$ and $Y_{o}\subseteq Y$ be Borel sets such that
$E(X_{o})=\overline{[\mathcal{U}^{\ast}\mathcal{H}_{2}]}$ and $F(Y_{o}%
)=\overline{[\mathcal{UH}_{1}]}$. Theorem \ref{210} shows that $\mathcal{S}%
_{1o}\mathcal{=S}_{1\chi}|_{E(X_{o})}$ and $\mathcal{S}_{2o}\mathcal{=S}%
_{2\chi}|_{F(Y_{o})}$ are complete orthocomplemented lattices, hence (since
they are commutative) complete Boolean lattices, and that $\chi$ induces a
complete Boolean lattice isomorphism $\chi_{o}$ between them. Denote by
$\mathfrak{A}_{i}$ the w*-closed algebra generated by $\mathcal{S}_{io}$
($i=1,2$), and let $\psi:\mathfrak{A}_{1}\rightarrow\mathfrak{A}_{2}$ be the
*-isomorphism induced by $\chi_{o}$. Now the $\mathfrak{A}_{i}$ are abelian
separably acting von Neumann algebras, hence there exist standard Borel
probability spaces $(X_{1},\mu_{1})$ and $(Y_{1},\nu_{1})$ such that the
$\mathfrak{A}_{i}$ are *-isomorphic to the corresponding $L^{\infty}$ spaces.
Moreover, we can take both $X_{1}$ and $Y_{1}$ to be compact intervals, which
we take to be $[\frac{1}{3},\frac{2}{3}]$ for notational convenience. Now the
inclusion $\mathfrak{A}_{1}\rightarrow\mathcal{D}_{1}$ maps the identity of
$\mathfrak{A}_{1}$ to $E(X_{o})$, and hence induces an injective unital
*-homomorphism $\theta_{1}:L^{\infty}(X_{1},\mu_{1})\rightarrow L^{\infty
}(X_{o},\mu)$. This map is implemented (see \cite{sic}) by a Borel function
$f_{1}:X_{o}\rightarrow X_{1}$; thus every projection $E\in\mathfrak{A}_{1}$ is of
the form $E=E(f_{1}^{-1}(\sigma))$ for some Borel subset $\sigma\subseteq
X_{1}$. Since $\theta_{1}$ is injective, $f_{1}$ can be taken to be onto
$X_{1}$. Similarly, there exists a Borel onto function $g_{1}:Y_{o}\rightarrow
Y_{1}$ such that every projection $F\in\mathfrak{A}_{2}$ is of the form
$F=F(g_{1}^{-1}(\tau))$ for some Borel subset $\tau\subseteq Y_{1}$. Now the
*-isomorphism $\psi:\mathfrak{A}_{1}\rightarrow\mathfrak{A}_{2}$ is also implemented
by a Borel bijection $h:Y_{1}\rightarrow X_{1}$; that is, for every projection
$E(f_{1}^{-1}(\sigma))$ we have $\psi(E(f_{1}^{-1}(\sigma)))=F(g_{1}%
^{-1}(h^{-1}(\sigma)))$. Define $f:X\rightarrow\lbrack0,1]$ by $f(x)=f_{1}(x)$
for $x\in X_{o}$ and $f(x)=0$ otherwise, and define $g:Y\rightarrow
\lbrack0,1]$ by $g(y)=h(g_{1}(y))$ for $y\in Y_{o}$ and $g(y)=1$ otherwise.
Thus $\mathcal{S}_{1\chi}=\{E(f^{-1}(\sigma)):\sigma\subseteq\lbrack0,1]$
Borel, $0\in\sigma\}$ and $\chi(E(f^{-1}(\sigma))=F(g^{-1}(\sigma))$.

Let $\kappa=\{(x,y)\in X\times Y:f(x)=g(y)\}$. Pick a countable dense subset
$\{s_{n}\}\subseteq\lbrack0,1)$ and note that $\kappa=\{(x,y)\in X\times
Y:\ $for each $n,\ x\in\alpha_{n}\Leftrightarrow y\in\beta_{n}\}$ where
$\alpha_{n}=f^{-1}([0,s_{n}])$ and $\beta_{n}=g^{-1}([0,s_{n}])$, so that
$\chi(E(\alpha_{n}))=F(\beta_{n})$ and $\chi(E(\alpha_{n}^{c}))\subseteq
F(\beta_{n}^{c})$. Thus, the complement of $\kappa$ can be written $\kappa
^{c}=\cup_{k}\gamma_{k}\times\delta_{k}$ where $\chi(E(\gamma_{k}))\subseteq
F(\delta_{k}^{c})$ for each $k$.

It is shown in the proof of Theorem 4.2 of \cite{eks} that an operator $T$ is
supported by $\kappa$ if and only if $F(\delta_{k})TE(\gamma_{k})=0$ for all
$k\in\mathbb{N}$. Thus if $T\in\mathcal{U}$ then $T$ is supported by $\kappa$.

Conversely, suppose $T$ is supported by $\kappa$. Then for each
$E=E(f^{-1}(\sigma))\in\mathcal{S}_{1\chi}$ we have
$\chi(E)^{\perp}TE=F(g^{-1}(\sigma)^{c})TE(f^{-1}(\sigma))$ since
$\chi(E)=F(g^{-1}(\sigma))$. But
the rectangle $f^{-1}(\sigma)\times g^{-1}(\sigma)^{c}$ is disjoint from
$\kappa$, hence $\chi(E)^{\perp}TE=0$ so that $T\in\mathcal{U}$.

\medskip

Finally, if $\mathcal{U}$ consists of all operators supported
by a set $\kappa$ of the above form, then writing
$\kappa=\{(x,y)\in X\times Y: \mbox{ for each }
n,\ x\in\alpha_{n}\Leftrightarrow y\in\beta_{n}\}$ as above
we have that  an operator $T$ is
in $\mathcal{U}$ if and only if $F(\beta_{n})^{\bot}TE(\alpha_{n})=0$ and
$F(\beta_{n})TE(\alpha_{n})^{\bot}=0$ for each $n$. Thus $\mathcal{U}%
=\{T\in\mathcal{B(H}_{1},\mathcal{H}_{2}):TE(\alpha_{n})=F(\beta_{n})T$ for
each $n\}$, and this is clearly closed under the triple product.
$\diamondsuit$\hfill

\begin{corollary}
In the notation of the last Theorem, if the algebra $[\mathcal{UU}^{\ast}]$
(resp. $[\mathcal{U}^{\ast}\mathcal{U}]$) is abelian, then the $\omega
$-support of $\mathcal{U}$ is the graph (resp. `reverse graph') of a Borel
function $f:X_{o}\rightarrow Y$ (resp. $g:Y_{o}\rightarrow X$)
for a suitable Borel subset $X_o \subseteq X$ (resp. $Y_o \subseteq Y$).
\end{corollary}

\noindent\ \textit{Proof }If $[\mathcal{U}^{\ast}\mathcal{U}]$ is abelian,
then its w*-closure is maximal abelian on $\overline{[\mathcal{U}^{\ast
}\mathcal{H}_{2}]}$. Hence it must be unitarily equivalent to the
multiplication algebra of $L^{\infty}(X_{o},\mu)$. It follows as in the above
proof (essentially by replacing $f_1$ by the identity  $X_o \to X_o$)
that there exists a Borel function $g:Y_{o}\rightarrow X$ such that the
$\omega$-support of $\mathcal{U}$ is $\{(x,y)\in X\times Y_{o}:x=g(y)\}$.

Similarly, if  $[\mathcal{UU}^{\ast}]$ is abelian, there exists a Borel
function $f:X_{o}\rightarrow Y$ such that the $\omega$-support of
$\mathcal{U}$ is $\{(x,y)\in X_{o}\times Y:f(x)=y\}$. $\diamondsuit$\hfill

\bigskip

\noindent\textbf{Remarks (i)} The previous Theorem is based on an idea
of V.S. Shulman. We take this opportunity to thank him.

\noindent\textbf{(ii)} Based on the characterization of Theorem \ref{th6} a
constructive proof of Theorem \ref{p1} can be given, avoiding the use
of the implementability of *-homomorphisms between $L^{\infty}$-spaces
(see \cite{it}).

\noindent\textbf{(iii)} In general, one function is not enough to describe the
support of a normalizing masa-bimodule. For an example, consider the von
Neumann algebra $\mathcal{A}=M_2(\m)$, where $\m$ is the
multiplication algebra of $L^{\infty}(0,1)$. It is not hard to see
that the support of $\mathcal{A}$ is a set of the form $\{(x,y):f(x)=f(y)\}$,
where $f$ is a certain Borel function on $[0,2]$ with period $1$.

\bigskip

It is easy to see that every finite rank operator in a von Neumann algebra
with abelian commutant is a sum of rank one operators in the algebra. Hence,
for normalizing masa bimodules Corollary \ref{prop88} can be improved as follows:

\begin{proposition}
\label{prop8} Let ${\mathcal{U}}$ be an ultraweakly closed normalizing masa
bimodule and suppose that $K$ is a compact operator in ${\mathcal{U}}$. Then
$K$ can be approximated in the norm topology by sums of rank one operators in
${\mathcal{U}}$. Moreover, if $K\in{\mathcal{C}}_{p}$ then it can be
approximated by sums of rank one operators in ${\mathcal{U}}$ in the $p$-norm
topology. Finally, every operator of rank $n$ in $\mathcal{U}$ is the sum of
$n$ rank one operators in $\mathcal{U}$.
\end{proposition}

Next we want to identify the normalizing masa-bimodules which are strongly
reflexive. If a rank one operator $R$ belongs to the masa-bimodule
${\mathcal{U}}$ and if $E_{R}$ and $F_{R}$ are the smallest projections in the
masas such that $F_{R}R=R$ and $RE_{R}=R$, then $F_{R}{\mathcal{B}%
}({{\mathcal{H}}_{1}},{{\mathcal{H}}_{2}})E_{R}\subseteq{\mathcal{U}}$. It was
proved in \cite{eks} that, if $\mathcal{U}$ is a strongly reflexive
masa-bimodule, the space
\[
\mbox{span }\{F_{R}{\mathcal{B}}({{\mathcal{H}}_{1}},{{\mathcal{H}}_{2}}%
)E_{R}:R\in{\mathcal{U}}\}
\]
is weakly dense in ${\mathcal{U}}$. Note that $\mathcal{B}({{\mathcal{H}}_{1}%
},{{\mathcal{H}}_{2}})$ is a strongly reflexive normalizing masa bimodule.
More generally, if $\{E_{n}\}\subseteq\mathcal{D}_{1}$ and $\{F_{n}%
\}\subseteq\mathcal{D}_{2}$ are countable families of mutually orthogonal
projections, it is easy to verify that the direct sum $\oplus_{n}%
F_{n}{\mathcal{B}}({{\mathcal{H}}_{1}},{{\mathcal{H}}_{2}})E_{n}$ is closed
under the triple product, and is clearly a strongly reflexive normalizing masa
bimodule. In fact, there are no others:

\begin{theorem}
\label{th11} A strongly reflexive subspace ${\mathcal{U}}$ is normalizing if
and only if there are countable families $\{E_{n}\}$ and $\{F_{n}\}$
consisting of mutually orthogonal projections, such that
\[
{\mathcal{U}}=\oplus F_{n}{\mathcal{B}}({{\mathcal{H}}_{1}},{{\mathcal{H}}%
_{2}})E_{n}.
\]
\end{theorem}

\noindent\textit{Proof. } Let $\mathcal{U}$ be a normalizing strongly
reflexive subspace and let $\u_o$ be its essential part acting on
$\k_1 = 0_+^{\bot}\hone$.
If $\mathcal{A}\subseteq\b(\k_{1})$ is
the von Neumann algebra $(\u_o^{\ast}\u_o)^{\prime\prime}$, we have shown
(Lemma \ref{lemma10}) that
${\mathcal{L}}=\mathop{\rm Lat}{\mathcal{A}}$ is a totally atomic
commutative Boolean
lattice. Let $\mathcal{E}=\{E_{n}:n=1,2,\ldots\}$ be the set of atoms of
$\mathcal{L}$ considered as projections in $\hone$
(necessarily countable since ${{\mathcal{H}}_{1}}$ is
separable). Since $\sum E_{n}=0_+^{\bot}$
in the strong operator topology, each
$T\in\mathcal{U}$ can be written $T=\sum TE_{n}$. Setting $F_{n}%
=[\mathcal{U}E_{n}]=\chi(E_{n})$ (where $\chi= \mathop{\rm Map} (\mathcal{U}%
)$), we see that the $F_{n}$ are orthogonal since $\chi$ is an ortho-map.
Since $TE_{n}=F_{n}TE_{n}$, it follows that $T=\sum F_{n}TE_{n}$. On the other
hand, for each $n$, each $T\in$ $F_{n}{\mathcal{B}}({{\mathcal{H}}_{1}%
},{{\mathcal{H}}_{2}})E_{n}$ satisfies $\chi(L)^{\bot}TL=0$ when $L$ is some
$E_{m}$, and hence when $L$ is in $\mathcal{L}$. Thus $F_{n}{\mathcal{B}%
}({{\mathcal{H}}_{1}},{{\mathcal{H}}_{2}})E_{n}\subseteq\mathcal{U}$ for each
$n$, and this completes the proof. $\diamondsuit$\hfill\bigskip

The proof of the above theorem actually gives more.

\begin{corollary}
\label{c12} The ultraweak closure of the rank one subspace of a normalizing
subspace $\mathcal{U}$ $\subseteq\mathcal{B(H},\mathcal{K)}$ is of the form
$\oplus_{n}{\mathcal{B}}({\mathcal{H}}_{n},{\mathcal{K}}_{n})$, for some
mutually orthogonal families of subspaces $\mathcal{H}_{n}\subseteq
\mathcal{H}$ and ${\mathcal{K}}_{n}\subseteq\mathcal{K}$.
\end{corollary}

\noindent\textit{Proof. } Obviously the rank one subspace of $\mathcal{U}$ is
also a normalizing space and hence its ultraweak closure is a reflexive
(Theorem \ref{th0}) normalizing masa bimodule. Thus the last theorem applies.
$\diamondsuit$\hfill\bigskip

According to \cite{eks}, a masa bimodule is strongly reflexive if and only if
its $\omega$-support $\kappa$ is marginally equivalent to a countable union of
Borel rectangles. In the terminology of \cite{eks}, $\kappa$ is the $\omega
$-closure of its $\omega$-interior. In this terminology, Corollary \ref{c12}
says that the $\omega$-interior of the $\omega$-support of a normalizing
masa-bimodule can be written as a countable union of Borel rectangles
$\alpha_{n}\times\beta_{n}$, such that the families $\{\alpha_{n}\}$ and
$\{\beta_{n}\}$ consist of disjoint Borel sets.

Let us recall that in \cite{eks} it was shown that although the rank one
subspace of a general strongly reflexive masa-bimodule is dense in the
bimodule in the weak operator topology, it need not be dense in the ultraweak
topology. By Theorem \ref{th11}, this cannot occur for normalizing masa-bimodules.

\begin{corollary}
If a normalizing subspace $\mathcal{U}$ is strongly reflexive, then its rank
one subspace is ultraweakly dense in $\mathcal{U}$.
\end{corollary}

\section{Normalizing masa-bimodules and synthesis}

Now we turn our attention to the question of spectral synthesis. Spectral
synthesis for operator algebras was introduced by Arveson \cite{a}. It can be
generalized for masa-bimodules as follows. If $\mathcal{L}_{1}=\mathcal{P(D}%
_{1}),\ \mathcal{L}_{2}=\mathcal{P(D}_{2})$ are the projection lattices of two
masas, a map $\varphi\in\mathcal{M(L}_{1},\mathcal{L}_{2})$ (a
\emph{commutative subspace map }in the terminology of \cite{opmap}) is said to
be \textbf{synthetic}, if the only ultraweakly closed masa-bimodule
${\mathcal{S}}$ with the property $\mathop{\rm
Map}{\mathcal{S}}=\varphi$ is $\mathop{\rm Op}\varphi$. There is a non trivial
fact hidden behind this definition. As was proved by Arveson (in fact, for the
case of CSL algebras, but it easily follows for masa-bimodules as well
\cite{dav}), given $\varphi\in\mathcal{M(L}_{1},\mathcal{L}_{2})$ there exists
an ultraweakly closed masa-bimodule ${\mathcal{M}}_{min}$, minimal with
respect to the property that its reflexive hull equals $\mathop{\rm Op}%
\varphi$. Thus a reflexive masa-bimodule ${\mathcal{M}}$ has synthetic map if
and only if ${\mathcal{M}}={\mathcal{M}}_{min}$. Von Neumann algebras with
abelian commutant have this property. The same holds for their generalization,
normalizing masa-bimodules.

\begin{theorem}
\label{th.14} Commutative subspace ortho-maps are synthetic.
\end{theorem}

This theorem will follow from a general fact about masa-bimodules (Proposition
\ref{prop16} below), which states that for every family of ultraweakly closed
masa-bimodules with the same reflexive hull, there is a certain natural
normalizing masa-bimodule contained in each member of the family. This
bimodule corresponds to the diagonal of a CSL algebra.

\begin{proposition}
\label{prop16} Let ${\mathcal{U}}=\mathop{\rm Op}\varphi\subseteq
{\mathcal{B}(\mathcal{H}_{1},\mathcal{H}_{2})}$ be a reflexive masa bimodule,
and let ${\mathcal{U}}_{0}=\mathop{\rm Op}\varphi\cap\mathop{\rm
Op}\varphi^{\bot}$. Then ${\mathcal{U}}_{0}$ is contained in ${\mathcal{U}%
}_{min}$, the minimal ultraweakly closed masa-bimodule with reflexive cover
${\mathcal{U}}$ .
\end{proposition}

\noindent\textit{Proof. } We consider the Hilbert space
$\mathcal{H=H}_2\oplus\mathcal{H}_1$ and the set of projections
\[
{\mathcal{S}}=\{\varphi(L)\oplus L:L\in\mathcal{S}_{1\varphi}\}
\]
on it. We put ${\mathcal{A}}= \mathop{\rm Alg} {\mathcal{S}_{1\varphi}}$,
${\mathcal{B}}= \mathop{\rm Alg} \varphi({\mathcal{S}_{1\varphi}})=
\mathop{\rm Alg} {\mathcal{S}_{2\varphi}}$ and $\mathcal{V=} \mathop{\rm Op}
\varphi^{\bot}$. \medskip

It is an easy verification to show that
$\Alg \mathcal{S} =
\left(  \smallmatrix \mathcal{B}  &  \mathcal{U} \\
\mathcal{V}^*  & \mathcal{A} \endsmallmatrix \right)$.
A similar calculation shows also that the diagonal $\mathcal{C\cap C}^{\ast}$
of $\mathcal{C}= \mathop{\rm Alg} \mathcal{S}$ (i.e. the commutant of
$\mathcal{S}$) is
${\mathcal{S}}^{\prime}=\left( \smallmatrix
{\mathcal{S}}_{2\varphi}^{\prime} & {\mathcal{U}}_{0}\\
{\mathcal{U}}_{0}^{\ast} & {\mathcal{S}}_{1\varphi}^{\prime}%
\endsmallmatrix \right)$.
But Arveson \cite{a} has shown that, for a CSL algebra $\mathcal{C}$, the
diagonal $\mathcal{C\cap C}^{\ast}$ is contained in $\mathcal{C}_{\min}$.
Therefore, to show that $\mathcal{U}_{0}\subseteq\mathcal{U}_{\min}$, it
suffices to prove that $\mathcal{C}_{\min}$ is contained in
\[
{\mathcal{M}}=\left(
\begin{array}
[c]{ll}%
{\mathcal{B}}_{\min} & {\mathcal{U}}_{\min}\\
{\mathcal{V}}_{\min}^{\ast} & {\mathcal{A}}_{\min}%
\end{array}
\right)  .
\]

Since $\mathcal{M}$ is a w*-closed subspace of $\mathcal{B}(H)$ containing the
masa $\mathcal{D}_{2}\oplus\mathcal{D}_{1}$, to show that
$\mathcal{C}_{\min}\subseteq\mathcal{M}$ it suffices to show that
$[\mathcal{M}\xi]=[\mathcal{C}\xi]$ for all $\xi\in \mathcal{H}$ (see
Theorem 22.19 of \cite{dav}). But this can be immediately verified using
the facts that $\mathcal{B}=\mathop{\rm Ref}(\mathcal{B}_{\min})$,
$\mathcal{U}=\mathop{\rm Ref}(\mathcal{U}_{\min})$,
$\mathcal{V}=\mathop{\rm Ref}(\mathcal{V}_{\min})$ and
$\mathcal{A=}\mathop{\rm Ref}(\mathcal{A}_{\min}).
\qquad \diamondsuit$\hfill\bigskip

\noindent\textit{Proof of Theorem \ref{th.14}.} Suppose that $\chi$ is a
commutative subspace map, which is also an ortho-map. The space $\mathcal{U}%
=\mathop{\rm Op} \chi$ is normalizing. But then $\mathcal{U}=\mathop{\rm Op}
\chi\cap\mathop{\rm Op} \chi^{\bot}$ (Theorem \ref{th0.1}), hence
$\mathcal{U}\subseteq\mathcal{U}_{min}$ by Proposition \ref{prop16}. Thus
$\mathcal{U=U}_{\min}$, which means that $\chi$ is synthetic. $\diamondsuit
$\hfill\bigskip

\begin{corollary}
\label{c16.1} If the map of an ultraweakly closed masa-bimodule $\mathcal{U}$
is an ortho-map, then ${\mathcal{U}}$ is normalizing.
\end{corollary}

\noindent\textit{{Proof. }} Since ${\mathcal{U}}$ is a masa-bimodule, its map
$\chi$ is commutative, hence by Theorem \ref{th.14} it is synthetic. Hence
$\mathcal{U=} \mathop{\rm Ref }\mathcal{U}=\mathop{\rm Op}\chi$. But since
$\chi$ is an ortho-map, $\mathop{\rm Op} \chi$ must be normalizing
(Theorem \ref{th0.1}). $\diamondsuit$\hfill\bigskip

We conclude this section with another result on synthesis. Let us call an
ultraweakly closed masa bimodule \emph{synthetic} if its map is synthetic.
Thus a synthetic masa-bimodule is automatically reflexive. As we will easily
see, normalizing masa-bimodules are not only synthetic, but in a sense
hereditarily synthetic.

\begin{proposition}
\label{16.1} \label{16.2} Let $\mathcal{U}$ be an ultraweakly closed
normalizing masa bimodule and let $\mathcal{S}\subseteq\mathcal{U}$ be an
ultraweakly closed masa-bimodule. If $\mathcal{S}$ is a
$\mathcal{U}^*\mathcal{U}$-submodule (or a $\mathcal{UU}^*$-submodule), then
$\mathcal{S}$ is synthetic.
\end{proposition}

\noindent\textit{Proof. } It suffices to show that ${\mathcal{S}}$ is
normalizing. Let $\mathcal{U}=\mathop{\rm Op} \varphi$, $\mathcal{B}_{1}
={\mathcal{U}}^{*}{\mathcal{U}}$ and $\mathcal{B}_{2}={\mathcal{U}%
}{\mathcal{U}}^{*}$. Suppose that $R,S,T\in\mathcal{S}$. Then $S,T\in
\mathcal{U}$ and hence $S^{\ast}T\in{\mathcal{B}}_{1}$. If $\mathcal{SB}%
_{1}\mathcal{\subseteq S}$, it follows that $RS^{\ast}T\in\mathcal{S}$. If
$\mathcal{B}_{2}\mathcal{S\subseteq S}$, we have $RS^{\ast}\in\mathcal{B}_{2}$
and hence $RS^{\ast}T\in\mathcal{S}$ again. $\diamondsuit$\hfill

\bigskip

An immediate application of the last proposition is the fact that, if
$\mathcal{A}$ is a von Neumann algebra with abelian commutant, every
ultraweakly closed left or right ideal of $\mathcal{A}$ is synthetic. Note
that this can also be inferred from the results of \cite{ls}.

\section{Normalizers of reflexive algebras}

In this section we present the relation between normalizing reflexive
subspaces and normalizers of reflexive algebras of operators. Let $\alga$
and $\b$ be reflexive algebras of operators on $\hone$ and $\htwo$
respectively. An operator $T\in \b(\hone,\htwo)$
is called a \textbf{semi-normalizer} of ${\mathcal{B}}$ into $\alga$ if
\begin{equation*}
T^{\ast }{\mathcal{B}}T\subseteq {\mathcal{A}}.
\end{equation*}
A \textbf{normalizer} of the algebra ${\mathcal{B}}$ into the algebra
$\mathcal{A}$ is a semi-normalizer of $\b$ into ${\mathcal{A}}$
whose adjoint is a semi-normalizer of ${\mathcal{A}}$ into ${\mathcal{B}}$.
We denote by $SN({\mathcal{B}},{\mathcal{A}})$ the set of semi-normalizers
of ${\mathcal{B}}$ into ${\mathcal{A}}$ and by $N(\b,\mathcal{A})$
the set of normalizers of ${\mathcal{B}}$ into ${\mathcal{A}}$.

In \cite{coates} it is shown that the set of normalizers of a nest algebra
into itself is closed in the weak operator topology. We show that, for
arbitrary reflexive algebras $\mathcal{A}$ and $\mathcal{B}$, the sets
$SN({\mathcal{B}},{\mathcal{A}})$ and $N({\mathcal{B}},{\mathcal{A}})$ are
closed in the strong operator topology, but not always in the weak operator
topology; nevertheless, Coates' result remains valid whenever the algebras
$\mathcal{A}$ and $\mathcal{B}$ are strongly reflexive.

\begin{proposition}
If $\mathcal{A}$ and $\mathcal{B}$ are reflexive algebras, the sets
$SN({\mathcal{B}},{\mathcal{A}})$ and $N({\mathcal{B}},{\mathcal{A}})$ are
strongly closed. If $\mathcal{B}$ is strongly reflexive (resp. $\mathcal{A}$
and $\mathcal{B}$ are strongly reflexive) then $SN({\mathcal{B}},{\mathcal{A}%
})$ (resp. $N({\mathcal{B}},{\mathcal{A}})$) is weakly closed. The sets
$SN(\mathbb{C}I,\mathbb{C}I)$ and $N(\mathbb{C}I,\mathbb{C}I)$ are not weakly closed.
\end{proposition}

\noindent\textit{Proof.} Suppose that $T_{\nu}\longrightarrow T$ strongly and
that $T_{\nu}\in SN({\mathcal{B}},{\mathcal{A}})$. Then,
for each operator $B\in{\mathcal{B}}$ and vectors $x$ and $y$,
$\langle BT_{\nu}x,T_{\nu}y\rangle\longrightarrow\langle BTx,Ty\rangle$. This
means, of course, that $T_{\nu}^{\ast}BT_{\nu}\longrightarrow T^{\ast}BT$ in
the weak operator topology. Since $T_{\nu}^{\ast}BT_{\nu}\in{\mathcal{A}}$ for
each $\nu$ and ${\mathcal{A}}$ is weakly closed, being reflexive, it follows
that $T^{\ast}BT\in{\mathcal{A}}$. The proof that $N(\b,\mathcal{A})$
is strongly closed is similar.

Now suppose that $T_{\nu}\longrightarrow T$ weakly and that $T_{\nu}\in
SN({\mathcal{B}},{\mathcal{A}})$. First note that for each \emph{finite rank}
operator $F\in\mathcal{B}$, $T_{\nu}^{\ast}FT_{\nu}\longrightarrow T^{\ast}FT$
weakly (indeed, writing $F=F_{1}^{\ast}F_{2}$ where $F_{1},F_{2}$ have finite
rank, we have $F_{i}T_{\nu}\longrightarrow F_{i}T$ strongly for $i=1,2$). Thus
$T^{\ast}FT\in\mathcal{A}$. Now if $B\in\mathcal{B}$ is arbitrary and
$x\in\mathcal{H}_{2}$, by strong reflexivity there is a net $(F_{i})$ of
finite rank operators in $\mathcal{B}$ such that $F_{i}Tx\rightarrow BTx$ and
so $T^{\ast}F_{i}Tx\rightarrow T^{\ast}BTx$. Since each $T^{\ast}F_{i}T$ is in
$\mathcal{A}$ and $\mathcal{A}$ is reflexive, it follows that $T^{\ast}%
BT\in\mathcal{A}$. This shows that $SN({\mathcal{B}},{\mathcal{A}})$ is weakly
closed. Therefore, if $\mathcal{A}$ is also strongly reflexive, the same holds
for $N(\mathcal{B},\mathcal{A})=SN({\mathcal{B}},{\mathcal{A}})\cap
(SN(\mathcal{A},\mathcal{B}))^{\ast}$.

The final statement is a consequence of the fact that the sets of isometries
and unitaries are not weakly closed. For a specific example, we construct a
sequence $T_{n}\in N(\mathbb{C}I,\mathbb{C}I)$ converging weakly to an
operator $T$ which is not in $SN(\mathbb{C}I,\mathbb{C}I)$. Let $\{e_{k}%
:k\in\mathbb{Z}\}\cup\{e\}$ be an orthonormal basis of $\mathcal{H}$ and let
$U$ be the bilateral shift defined on $\mathcal{H}_{o}=\overline{[e_{k}%
:k\in\mathbb{Z}]}$ by $Ue_{k}=e_{k+1}$. Define $T_{n}=U^{n}\oplus
I\in\mathcal{B(H)}$ and $T=0 \oplus I \in \b(\h)$. Each
$T_{n}$ is a unitary operator on $\mathcal{H}$, and so $T_{n}^{\ast
}(\mathbb{C}I)T_{n}\subseteq\mathbb{C}I$ and $T_{n}(\mathbb{C}I)T_{n}^{\ast
}\subseteq\mathbb{C}I$, i.e. $T_{n}\in N(\mathbb{C}I,\mathbb{C}I)$. Since
$U^{n}\rightarrow0$ weakly, $T_{n}\rightarrow T$ weakly. But $T^{\ast
}(\mathbb{C}I)T=[T]\nsubseteqq\mathbb{C}I.
\qquad \diamondsuit$\hfill\bigskip

Let $\mathcal{L}_1=\Lat\mathcal{A}$ and
$\mathcal{L}_2=\Lat\mathcal{B}$.
If $\varphi \in {\mathcal{M}({\mathcal{L}_{1}},{\mathcal{L}_{2}})}$, we put
\begin{equation*}
{\mathcal{U}}_{\varphi }=\{T\in \mathcal{B}({{\mathcal{H}}_{1}},{{\mathcal{H}%
}_{2}}):TL=\varphi (L)T\mbox{
for all }L\in \mathcal{L}_{1}\}
\end{equation*}
which is a reflexive normalizing subspace, and also a bimodule over the
diagonals $\mathcal{B\cap B}^{\ast }$ and $\mathcal{A\cap A}^{\ast }$, as is
readily verified.

\begin{theorem}
\label{th17} The set $SN({\mathcal{B}},{\mathcal{A}})$ of semi-normalizers
of ${\mathcal{B}}$ into ${\mathcal{A}}$ is a union of reflexive normalizing
subspaces. More precisely,
\begin{equation*}
SN({\mathcal{B}},{\mathcal{A}})=\cup \{{\mathcal{U}}_{\varphi }:\varphi \in {%
\mathcal{M}({\mathcal{L}_{1}},{\mathcal{L}_{2}})}\}.
\end{equation*}
\end{theorem}

\noindent \textit{Proof. } Suppose that $\varphi \in {\mathcal{M}({\mathcal{L%
}_{1}},{\mathcal{L}_{2}})}$ and $T\in {\mathcal{U}}_{\varphi }$. Then, for
each $L\in {\mathcal{L}_{1}}$ and $B\in {\mathcal{B}}$,
\begin{equation*}
L^{\perp }T^{\ast }BTL=L^{\perp }T^{\ast }B\varphi (L)T=L^{\perp }T^{\ast
}\varphi (L)B\varphi (L)T=L^{\perp }LT^{\ast }B\varphi (L)T=0
\end{equation*}
since $\varphi (L)\in \mathcal{L}_{2}$, so that $T^{\ast }BT\in \mathcal{A}$%
, which shows that $T$ is a semi-normalizer of ${\mathcal{B}}$ into ${%
\mathcal{A}}$. Thus $\mathcal{U}_{\varphi }\subseteq SN(\mathcal{B},\mathcal{%
A})$.

Conversely, let $T\in SN(\mathcal{B},\mathcal{A})$. Define $\varphi \in
\mathcal{M(L}_{1},\mathcal{L}_{2})$ by
\begin{equation*}
\varphi (L)=\overline{[\mathcal{B}T(L)]},\ \ \ L\in \mathcal{L}_{1}.
\end{equation*}
It is obvious that $\varphi (L)\in \mathcal{L}_2$. It is also easy to see
that $\varphi $ is 0-preserving and~join-continuous. We will show that
$T\in {\mathcal{U}}_{\varphi }$. Let $L\in \mathcal{L}_1$. Since
$T(L{{\mathcal{H}}_{1}})\subseteq \mathcal{B}T(L{{\mathcal{H}}_{1}})$,
it is clear that $\varphi (L)^{\bot }TL=0$. On the other hand, for each
$B\in \mathcal{B}$ we have $T^{\ast }BT\in \mathcal{A}$ and so
$L^{\bot }T^{\ast }BTL=0$; hence $\langle BTLx,TL^{\bot }y\rangle =0$ for
each $x,y$. Since the closure of
$\{ BTLx : B \in \mathcal{B}, x \in \mathcal{H}_1 \}$ is
$\varphi (L)(\mathcal{H}_2)$, it follows that
$\varphi (L)(\mathcal{H}_2) \bot TL^{\bot}(\mathcal{H}_1)$.
Thus $\langle \varphi (L)z,TL^{\bot }y\rangle =0$ for all
$y\in \mathcal{H}_1, z\in \mathcal{H}_2$ and so
$\varphi (L)TL^{\bot }=0$. Adding the relation $\varphi (L)TL=TL $,
we obtain $\varphi (L)T=TL$ so that $T\in \mathcal{U}_{\varphi }$
as required. $\diamondsuit $\hfill \bigskip

It is obvious that the intersection of several normalizing reflexive
subspaces is again a normalizing reflexive subspace. If, for each $\psi\in{%
\mathcal{M}}({\mathcal{L}_{2}},{\mathcal{L}_{1}})$, we put
\begin{equation*}
{\mathcal{V}}_{\psi}=\{S\in{\mathcal{B}}({{\mathcal{H}}_{2}},{{\mathcal{H}}%
_{1}}):SM=\psi(M)S\mbox{ for all }\ M\in{\mathcal{L}_{2}}\},
\end{equation*}
the above theorem gives us the following

\begin{corollary}
\label{c18} The set $N({\mathcal{B}},{\mathcal{A}})$ of normalizers of ${%
\mathcal{B}}$ into ${\mathcal{A}}$ is a union of reflexive normalizing
subspaces. More precisely,
\begin{equation*}
N({\mathcal{B}},{\mathcal{A}})=\cup \{{\mathcal{U}}_{\varphi }\cap {\mathcal{%
V}}_{\psi }^{\ast }:\varphi \in {\mathcal{M}({\mathcal{L}_{1}},{\mathcal{L}%
_{2}})},\psi \in {\mathcal{M}}({\mathcal{L}_{2}},{\mathcal{L}_{1}})\}.
\end{equation*}
\end{corollary}

This corollary was proved in \cite{coates} for the case of nest algebras.

\medskip

Recall (Proposition \ref{prop-1}) that w*-closed normalizing subspaces
are generated in norm by their partial isometries. Thus, the previous
results have the following consequence:
\begin{corollary}
Any seminormalizer (resp. normalizer) of $\mathcal{B}$ into $\mathcal{A}$
is the norm-limit of linear combinations of partial isometries in
$SN(\b,\mathcal{A})$ (resp. in $N(\b,\mathcal{A})$).
\end{corollary}

Having shown that the set of normalizers between reflexive algebras is a
union of normalizing spaces, we now turn to the converse question: when does a
normalizing reflexive space ${\mathcal{U}}$ consist of semi-normalizers or
normalizers of two reflexive algebras ${\mathcal{A}}$ and ${\mathcal{B}}$?
For example, it is easily verified that every operator in ${\mathcal{U}}$ is
a normalizer of the algebra $\mathop{\rm Alg}{\mathcal{S}}_{2\chi}$ into the
algebra $\mathop{\rm Alg}{\mathcal{S}}_{1\chi}$ (where
$\chi=\mathop{\rm Map}\mathcal{U}$).

We write $\mathcal{A}_{d}=\mathcal{A\cap A}^{\ast}$and $\mathcal{B}_{d}=%
\mathcal{B\cap B}^{\ast}$. If $\mathcal{K}_{1}=\mathcal{U}^{\ast }(\mathcal{H%
}_{2})=0_{+}^{\bot}\mathcal{H}_{1}$ and $\mathcal{K}_{2}=\mathcal{U(H}%
_{1})=I_{-}\mathcal{H}_{2}$, we denote by $\mathcal{A}_{o}$ (resp. $\mathcal{%
B}_{o}$) the compression of $\mathcal{A}$ (resp. $\mathcal{B}$) to $\mathcal{%
K}_{1}$ (resp. $\mathcal{K}_{2}$); we write $\mathcal{U}_{o}\subseteq%
\mathcal{B(K}_{1},\mathcal{K}_{2})$ for the restriction of $\mathcal{U}$ to $%
\mathcal{K}_{1}$ and $\chi_{o}$ for the map of $\mathcal{U}_{o}$.

\begin{lemma}
\label{l54}The following are equivalent:

(1) $\mathcal{U}^{\ast }\mathcal{U\subseteq A}$.

(2) For each $L\in \mathop{\rm Lat}\mathcal{A}_{d}$, $\chi (L)\bot \chi
(L^{\bot })$.

(3) For each $L\in \mathop{\rm Lat}\mathcal{A}$, $\chi (L)\bot \chi (L^{\bot
})$.

(4) With respect to the decomposition $\mathcal{H}_{1}=\mathcal{K}_{1}\oplus
\mathcal{K}_{1}^{\bot }$, each $L\in \mathop{\rm Lat}\mathcal{A}$ decomposes
as $L=L_{1}\oplus L_{2}$, where $L_{1}\in \mathcal{S}_{1\chi _{o}}$.
\end{lemma}

\proof Since $\mathcal{U}^{\ast }\mathcal{U\subseteq A}$ is equivalent to $%
\mathcal{U}^{\ast }\mathcal{U\subseteq A}_{d}$ and hence to $\mathop{\rm Lat}%
\mathcal{A}_{d}\mathcal{\subseteq }\mathop{\rm Lat}\mathcal{U}^{\ast }%
\mathcal{U}$ and to $\mathop{\rm Lat}\mathcal{A\subseteq }\mathop{\rm Lat}%
\mathcal{U}^{\ast }\mathcal{U}$, the Lemma is an immediate application of
Theorem \ref{210}, since $\mathop{\rm Lat}{\mathcal{U}}^{\ast }{\mathcal{U}}%
=\{L\in {\mathcal{P}}_{1}:\chi (L)\bot \chi (L^{\bot })\}=\{L_{1}\oplus
L_{2}:L_{1}\in {\mathcal{S}}_{1\chi _{o}}\}.$ $\diamondsuit $\hfill \bigskip

\begin{theorem}
\label{char_sn} Let $\mathcal{U\subseteq B(H}_{1},\mathcal{H}_{2})$ be a
w*-closed normalizing space. The following are equivalent:

(1) $\mathcal{U}\subseteq SN(\mathcal{B},\mathcal{A)}$.

(2) (i) $\chi (L)\bot \chi (L^{\bot })$ for each $L\in \Lat \mathcal{A}$
and (ii) $\chi _{o}( \Lat \mathcal{A}_{o})\subseteq \Lat \mathcal{B}_{o}$.

(3) (i) $\chi (L)\bot \chi (L^{\bot })$ for each $L\in \Lat \mathcal{A}$
and (ii) for each $L\in \mathop{\rm Lat}\mathcal{A}$, there
exists a projection $Q$ such that $\overline{[\mathcal{B}\chi (L)]}=\chi
(L)|_{\mathcal{K}_{1}}\oplus Q$.
\end{theorem}

\proof\textbf{(1)}$\Leftrightarrow$\textbf{(3)} Suppose (1) holds. Then $%
\mathcal{U}^{\ast}\mathcal{U\subseteq A}$, and (i) holds by the Lemma. For
all $T,S\in\mathcal{U}$, $B\in\mathcal{B}$  and $L \in \Lat \alga$ we have
\begin{equation*}
\langle BSLx,TL^{\bot}y\rangle=0\qquad(x,y\in\mathcal{H}_{1}).
\end{equation*}
Thus $\overline{[\mathcal{B}\chi(L)]}\leq\chi(L^{\bot})^{\bot}=(\chi(L)^{%
\bot }\wedge\chi(I))^{\bot}=\chi(L)\vee\chi(I)^{\bot}$. Writing $%
L=L_{1}\oplus L_{2}$ where $L_{1}\in\mathcal{S}_{1\chi_{o}}$ we have $%
\chi(L)=\chi_{o}(L_{1})\oplus0$. Since $\chi(L)\leq\overline{\lbrack\mathcal{%
B}\chi(L)]}$ we obtain
\begin{equation*}
\chi_{o}(L_{1})\oplus0\leq\overline{\lbrack\mathcal{B}\chi(L)]}\leq\chi
_{o}(L_{1})\oplus I
\end{equation*}
which shows that $\overline{[\mathcal{B}\chi(L)]}=\chi_{o}(L_{1})\oplus Q$,
for some projection $Q$, as required.

These steps can clearly be reversed.

\medskip

\textbf{(1)}$\Rightarrow $\textbf{(2)} Again (i) holds by the Lemma. For
(ii), note first that $\mathcal{U}_{o}\subseteq SN(\mathcal{B}_{o},\mathcal{A%
}_{o})$. Thus for each $L\in \mathop{\rm Lat}\mathcal{A}_{o}$, $S,T\in
\mathcal{U}_{o}$ and $B\in \mathcal{B}_{o}$ we have $L^{\bot }T^{\ast }BSL=0$%
. As in the proof of (1)$\Rightarrow $(3), we conclude that $\chi
_{o}(L)\leq \overline{[\mathcal{B}_{o}\chi _{o}(L)]}\leq \chi _{o}(L^{\bot
})^{\bot }=\chi _{o}(L)$ since $\chi _{o}(I)=I$. Thus $\chi _{o}(L)=%
\overline{[\mathcal{B}_{o}\chi _{o}(L)]}\in \mathop{\rm Lat}\mathcal{B}_{o}$.

\medskip

\textbf{(2)}$\Rightarrow$\textbf{(1)} Given $B\in\mathcal{B}$ and $T\in%
\mathcal{U}$, we will prove that $T^{\ast}BT\in\mathcal{A}$, equivalently
that $L^{\bot}T^{\ast}BTL=0$ for each $L\in\mathop{\rm Lat}\mathcal{A}$.
Since $\mathcal{U}^{\ast}\mathcal{U\subseteq A}$, we have $L=L_{1}\oplus
L_{2}$ where $L_{1}\in\mathcal{S}_{1\chi_{o}}$ by the Lemma. It is easily
seen that $L_{1}\in\mathop{\rm Lat}\mathcal{A}_{o}$.

Since $T\in\mathcal{U}$, we have $T(L^{\bot}\mathcal{H}_{1})\subseteq
\chi(L^{\bot})\mathcal{H}_{2}\subseteq\chi(L)^{\bot}\mathcal{H}_{2}$ since $%
\chi(L^{\bot})\bot\chi(L)$ by the Lemma. Thus $\chi(L)TL^{\bot}=0$. Since
also $\chi(L)^{\bot}TL=0$, it follows that $TL=\chi(L)T$ and $%
TL^{\bot}=\chi(L)^{\bot}\chi(I)T$ since the range of $T$ is contained in $%
\chi(I)$. Thus $L^{\bot}T^{\ast}BTL=T^{\ast}\chi(L)^{\bot}\chi(I)B\chi(I)%
\chi(L)T$. But $\chi(I)B|_{\mathcal{K}_{2}}\in\mathcal{B}_{o}$ and since $%
\chi (L)|_{\mathcal{K}_{2}}=\chi_{o}(L_{1})\in\mathop{\rm
Lat}\mathcal{B}_{o}$ by assumption we have $(\chi(I)B\chi(I))\chi
(L)=\chi(L)(\chi(I)B\chi(I))\chi(L)$ and so $L^{\bot}T^{\ast}BTL=0$. This
shows that $T^{\ast}BT\in\mathcal{A}$ and concludes the proof. $\diamondsuit$

\begin{corollary}
\label{cor_ess_sn} Let ${\mathcal{U}}\subseteq \mathcal{B(H}_{1},\mathcal{H}%
_{2})$ be an essential normalizing ultraweakly closed subspace and $\chi =%
\mathop{\rm
Map}{\mathcal{U}}$. Then ${\mathcal{U}}\subseteq SN({\mathcal{B}},{\mathcal{A%
}})$ if and only if

\emph{(i)} $\chi (L)\bot \chi (L^{\bot })$ for each $L\in \Lat \alga$ and

\emph{(ii)} $\chi (\mathop{\rm Lat}\mathcal{A})\subseteq \mathop{\rm Lat}%
\mathcal{B}$.
\end{corollary}

\noindent\textit{Proof. } Immediate from Theorem \ref{char_sn}. $%
\diamondsuit $

\begin{corollary}
\label{cor_n} Let ${\mathcal{U}}\subseteq \mathcal{B(H}_{1},\mathcal{H}_{2})$
be a normalizing ultraweakly closed subspace and let
$\chi =\mathop{\rm Map}{\mathcal{U}}$ and
$\chi ^{\ast }=\mathop{\rm Map}{\mathcal{U}}^{\ast }$.
Then ${\mathcal{U}}\subseteq N({\mathcal{B}},{\mathcal{A}})$ if and only if

\emph{(i)} $\chi (L)\bot \chi (L^{\bot })$ for each $L\in \Lat \alga$
while $\chi ^{\ast }(M)\bot \chi ^{\ast }(M^{\bot })$ for
each $M\in \mathop{\rm Lat}\mathcal{B}$ and

\emph{(ii)} $\chi _{o}(\mathop{\rm Lat}\mathcal{A}_{o})=\mathop{\rm Lat}%
\mathcal{B}_{o}$.
\end{corollary}

\noindent \textit{Proof. } If ${\mathcal{U}}\subseteq N({\mathcal{B}},{%
\mathcal{A}})$, then $\chi _{o}(\mathop{\rm Lat}\mathcal{A}_{o})\subseteq %
\mathop{\rm Lat}\mathcal{B}_{o}$ and $\chi _{o}^{\ast }(\mathop{\rm Lat}%
\mathcal{B}_{o})\subseteq \mathop{\rm Lat}\mathcal{A}_{o}$ by Theorem \ref
{char_sn}. But ${\mathcal{U}}_{o}^{\ast }{\mathcal{U}}_{o}\subseteq {%
\mathcal{A}}_{o}$ and hence $\mathop{\rm Lat}{\mathcal{A}}_{o}\subseteq %
\mathop{\rm Lat}({\mathcal{U}}_{o}^{\ast }{\mathcal{U}}_{o})={\mathcal{S}}%
_{1\chi _{o}}$ (Theorem \ref{210}). It follows that $\chi _{o}$ is
one-to-one on $\mathop{\rm Lat}{\mathcal{A}}_{o}$ into $%
\mathop{\rm
Lat}{\mathcal{B}}_{o}$. Similarly, $\mathop{\rm Lat}{\mathcal{B}}%
_{o}\subseteq \mathop{\rm
Lat}({\mathcal{U}}_{o}{\mathcal{U}}_{o}^{\ast })={\mathcal{S}}_{2\chi _{o}}$
and hence the inverse of $\chi _{o}$, namely $\chi _{o}^{\ast }$, maps $%
\mathop{\rm
Lat}{\mathcal{B}}_{o}$ 1-1 into $\mathop{\rm Lat}{\mathcal{A}}_{o}$.

Suppose conversely that \emph{(i)} and \emph{(ii)} hold. By Theorem \ref
{char_sn}, $\mathcal{U\subseteq }SN(\mathcal{B},\mathcal{A})$. As above,
\emph{(i) }gives $\mathop{\rm
Lat} {\mathcal{A}}_{o}\subseteq {\mathcal{S}}_{1\chi _{o}}$ and hence $\chi
_{o}^{\ast }\chi _{o}|_{\mathop{\rm Lat}{\mathcal{A}}_{o}}=id|_{%
\mathop{\rm
Lat}{\mathcal{A}}_{o}}$ from Theorem \ref{th0.001}. Therefore $\chi ^{\ast }(%
\mathop{\rm Lat}{\mathcal{B}}_{o})\subseteq \mathop{\rm Lat}{\mathcal{A}}%
_{o} $ and so ${\mathcal{U}}^{\ast }\subseteq SN({\mathcal{A}},{\mathcal{B}}%
) $ again from Theorem \ref{char_sn}. $\diamondsuit $\hfill \bigskip

\begin{corollary}
\label{char_ess_n} Let ${\mathcal{U\subseteq B(H}}_{1},\mathcal{H}_{2})$ be
an essential normalizing ultraweakly closed subspace and let $\chi =%
\mathop{\rm
Map}{\mathcal{U}}$. Then ${\mathcal{U}}\subseteq N({\mathcal{B}},{\mathcal{A}%
})$ if and only if

\emph{(i)} $\chi (L)\bot \chi (L^{\bot })$ for each $L\in \Lat \alga$ and

\emph{(ii)} $\chi (\mathop{\rm Lat}\mathcal{A})=\mathop{\rm Lat}\mathcal{B}$.
\end{corollary}

\noindent \textit{Proof. } From Corollary \ref{cor_n} it is sufficient to
show that $(i)$ and $(ii)$ imply that $\chi ^{\ast }(M)\bot \chi ^{\ast
}(M^{\bot })$ for each $M\in \mathop{\rm Lat}\mathcal{B}_{d}$. But this is
immediate from Lemma \ref{l54} (applied to $\mathcal{U}^{\ast }$) since from
$(ii)$ we have $\mathop{\rm Lat}\mathcal{B}=\chi (\mathop{\rm Lat}\mathcal{A}%
)\subseteq {\mathcal{S}}_{2}=\mathop{\rm Lat}(\mathcal{UU}^{\ast })$ and so $%
\mathcal{UU}^{\ast }\subseteq \mathcal{B}$. $\diamondsuit $\hfill \bigskip

Note that Theorem \ref{th17} yields, for every $T\in SN({\mathcal{B}},{%
\mathcal{A}})$, a certain normalizing subspace ${\mathcal{U}}_{\varphi }$,
such that $T\in {\mathcal{U}}_{\varphi }\subseteq SN({\mathcal{B}},{\mathcal{%
A}})$. We show that this property extends from single operators to \emph{%
linear spaces} consisting of seminormalizers. Recall that ${\mathcal{U}}%
_{\varphi }$ is a ${\mathcal{B}}_{d},{\mathcal{A}}_{d}$-bimodule, where ${%
\mathcal{A}}_{d}$ and ${\mathcal{B}}_{d}$ are the diagonals of ${\mathcal{A}}
$ and ${\mathcal{B}}$. It may be interesting to note that all right ${%
\mathcal{A}}_{d}$-modules $\mathcal{U}$ which consist of semi-normalizers of
${\mathcal{B}}$ into ${\mathcal{A}}$ automatically satisfy a rather strong
condition.

\begin{proposition}
\label{prop100}\textbf{(1)} Any \emph{linear space} $\mathcal{U}\subseteq SN(%
\mathcal{B},\mathcal{A)}$ (resp. $\mathcal{U}\subseteq N(\mathcal{B},%
\mathcal{A)}$) is contained in an ultraweakly closed normalizing space $%
\mathcal{U}_{A}\subseteq SN(\mathcal{B},\mathcal{A)}$ (resp.
$\u_A\subseteq N(\mathcal{B},\mathcal{A})$) which is a right
$\mathcal{A}_{d}$-module.

\textbf{(2) }If $\mathcal{U}\subseteq SN({\mathcal{B}},{\mathcal{A}})$ is a
right $\mathcal{A}_{d}$-module, then $\mathcal{U}$ is automatically
normalizing. Moreover, $0_{+}=\ker \mathcal{U}$ is in the centre of $%
\mathcal{A}_{d}$ and the restriction $\mathcal{A}_{od}$ of $\mathcal{A}_{d}$
to $\mathcal{K}_{1}=0_{+}\mathcal{H}_{1}$ is the bicommutant of $(\mathcal{U}%
^{\ast }\mathcal{U})|_{\mathcal{K}_{1}}$. It follows that $\mathcal{S}%
_{1\chi }=\{L\oplus I:L\in \mathop{\rm
Lat}\mathcal{A}_{od}\}$.
\end{proposition}

\proof \textbf{(1)} Suppose $\u \subseteq SN(\b, \alga)$ is a linear space.
Define $\u_a=[\mathcal{UA}_{d}]$ and
$\u_A=\overline{[\mathcal{UA}_{d}]}^{uw}$. If $T_{1},T_{2}\in \u$
and $A_{1},A_{2}\in \mathcal{A}_{d}$, then for each $B\in \b$
the fact that $T_i^*BT_j \in \mathcal{A}$ and $A_i,A_i^*\in \mathcal{A}$ yields
$(T_1A_1+T_2A_2)^*B(T_{1}A_{1}+T_{2}A_{2})\in \mathcal{A}$.
Thus $\mathcal{U}_{a}\subseteq SN(\mathcal{B},\mathcal{A})$.

Now $\mathcal{U}_{a}\mathcal{A}_{d}\subseteq \mathcal{U}_{a}$ by
construction and $\mathcal{U}_{a}^{\ast }\mathcal{U}_{a}\subseteq \mathcal{A}
$ since $\mathcal{U}_{a}\subseteq SN(\mathcal{B},\mathcal{A})$. Therefore $%
\mathcal{U}_{a}$ is normalizing (Remark \ref{rem1}).

It follows from Theorem \ref{24} that $\mathcal{U}_{A}$ is a reflexive
normalizing space. Also, since $\mathcal{U}_{A}$ is the strong operator
closure of $\mathcal{U}_{a}$ and $SN(\mathcal{B},\mathcal{A)}$ is strongly
closed, $\mathcal{U}_{A}\subseteq SN(\mathcal{B},\mathcal{A)}$. That the
strong closure of a right $\mathcal{A}_{d}$-module is a right $\mathcal{A}%
_{d}$-module is obvious.

The case $\mathcal{U}\subseteq N(\mathcal{B},\mathcal{A)}$
is similar.

\smallskip

\textbf{(2)} That $\mathcal{U}$ is normalizing follows from Remark \ref{rem1}
as above. Let $\chi =\mathop{\rm Map}\mathcal{U}$. Since $\mathcal{UA}_{d}%
\mathcal{\subseteq U}$ we have $\mathcal{S}_{1\chi }\subseteq
\mathop{\rm
Lat}\mathcal{A}_{d}$ (Lemma \ref{le1}). On the other hand since ${\mathcal{U}%
}\subseteq SN({\mathcal{B}},{\mathcal{A}})$ we have $\mathcal{U}^{\ast }%
\mathcal{U}\subseteq \mathcal{A}_{d}$. Thus $0_{+}$ commutes with $\mathcal{A%
}_{d}$ and also $0_{+}\in (\mathcal{U}^{\ast }\mathcal{U})^{\prime \prime }%
\mathcal{\subseteq A}_{d}$, so $0_{+}$ is in the centre of $\mathcal{A}_{d}$%
. Hence $\mathcal{A}_{d}$ can be written $\mathcal{A}_{od}\oplus \mathcal{A}%
_{1d}$ with respect to the decomposition $\mathcal{H}_{1}=\mathcal{K}%
_{1}\oplus \mathcal{K}_{1}^{\bot }$. Using Theorem \ref{210}, we have
\begin{equation*}
(\mathcal{U}^{\ast }\mathcal{U})^{\prime \prime }=\left(
\begin{array}{cc}
(\mathcal{S}_{1\chi }|_{\mathcal{K}_{1}})^{\prime } & 0 \\
0 & \mathbb{C}I
\end{array}
\right) \subseteq \mathcal{A}_{d}\subseteq \mathcal{S}_{1\chi }^{\prime
}=\left(
\begin{array}{cc}
(\mathcal{S}_{1\chi }|_{\mathcal{K}_{1}})^{\prime } & 0 \\
0 & \mathcal{B(K}_{1}^{\bot })
\end{array}
\right)
\end{equation*}
This shows that $\alga_{od}$ equals
$(\s_{1\chi }|_{\k_{1}})^{\prime }=(\u^*\u|_{\k_1})^{\prime \prime }$ and
$\Lat \alga_{od}$ equals $\s_{1\chi }|_{\k_{1}}= \Lat (\u^*\u|_{\k_{1}})$.
Thus $\s_{1\chi } =\{L\oplus I:L\in \Lat(\u^*\u|_{\k_{1}})\}
=\{L\oplus I:L\in \Lat\mathcal{A}_{od}\}.
\qquad \diamondsuit $\hfill \bigskip

It is clear that every linear space $\mathcal{U}\subseteq SN({\mathcal{B}},{%
\mathcal{A}})$ is also contained in a w*-closed normalizing space which is a
$(\mathcal{B}_{d},\mathcal{A}_{d})$-bimodule (just consider the w*-closure
of $[\mathcal{B}_{d}\mathcal{UA}_{d}]$). This bimodule is not necessarily
maximal with respect to being a linear space of seminormalizers.

For an example, consider the nest algebra $\mathcal{A}$ of all upper
triangular $2\times 2$ matrices and let $\mathcal{U}$ be the space of all
strictly lower triangular matrices. One easily checks that ${\mathcal{U}}%
\subseteq SN({\mathcal{A}},\mathcal{A})$ and that ${\mathcal{U}}$ is a
bimodule over the diagonal algebra (which is a masa in this case). However,
it is readily verified that the linear span of the matrix unit $E_{11}$ and $%
\mathcal{U}$ is also contained in $SN({\mathcal{A}},\mathcal{A})$.

\begin{remark}
Let ${\mathcal{U}}\subseteq SN({\mathcal{B}},{\mathcal{A}})$ be an essential
w*-closed normalizing space, which is a right ${\mathcal{A}}_{d}$-module. Then $\u$
is maximal with respect to being a linear space in $SN(\b,\mathcal{A}).$
\end{remark}

\proof Let $\mathcal{V}\subseteq SN({\mathcal{B}},{\mathcal{A}})$ be a
linear space containing $\mathcal{U}$ and let $T\in\mathcal{V}$. For each
$S\in\mathcal{U}$, we have
$S^{\ast}T\in\mathcal{V}^{\ast}\mathcal{V\subseteq A}_{d}$.
Since $\mathcal{UA}_{d}\mathcal{\subseteq U}$, if
$\chi=\mathop{\rm Map}{\mathcal{U}}$ we have
$\mathcal{S}_{1\chi}\subseteq\mathop{\rm Lat}\mathcal{A}_{d}$
(Lemma \ref{le1}). Thus $L^{\bot}S^{\ast}TL=0$ for each
$L\in\mathcal{S}_{1\chi}$ and hence $\langle TLx,SL^{\bot}y\rangle=0$ for
all $x,y$. 
But the closure of
$\{SL^{\bot}y:S\in\mathcal{U}, y \in \mathcal{H}_1 \}$ is
$\chi(L^{\bot}) \mathcal{H}_2$
which equals $\chi(L)^{\bot} \mathcal{H}_2$ since $\chi$ is
essential; therefore 
$TL(\mathcal{H}_1) \bot \chi(L)^{\bot} \mathcal{H}_2$. We have shown
that $\chi(L)^{\bot}TL=0$ for all $L\in\mathcal{S}_{1\chi}$ and so $T\in
\mathcal{U}$. $\diamondsuit$\hfill

\bigskip

Suppose now that ${\mathcal{A}}$ and ${\mathcal{B}}$ are CSL algebras, that
is, the respective invariant subspace lattices are commutative. In this case
Proposition \ref{prop8} and Theorem \ref{th.14} immediately yield the next
corollary. The last statement was proved by Coates \cite{coates} for the
case of nest algebras.

\begin{corollary}
\label{c19} If $\mathcal{A}$ and $\mathcal{B}$ are CSL algebras, the set $SN(%
{\mathcal{B}},{\mathcal{A}})$ ($N({\mathcal{B}},{\mathcal{A}})$) of
semi-normalizers (normalizers) of ${\mathcal{B}}$ \ into ${\mathcal{A}}$ is
a union of synthetic normalizing masa-bimodules. Each compact operator $K$
in $SN({\mathcal{B}},{\mathcal{A}})$ ($N({\mathcal{B}},{\mathcal{A}})$) can
be approximated in norm by sums of rank one operators in $SN({\mathcal{B}},{%
\mathcal{A}})$ ($N({\mathcal{B}},{\mathcal{A}})$). Moreover, if $K\in {%
\mathcal{C}}_{p}$, then it can be approximated in the ${\mathcal{C}}_{p}$%
-norm by sums of rank one operators in $SN({\mathcal{B}},{\mathcal{A}})$ ($N(%
{\mathcal{B}},{\mathcal{A}})$). Finally, if $K$ has finite rank, say $n$, it
can be written as a sum of $n$ rank one operators in $SN({\mathcal{B}},{%
\mathcal{A}})$ ($N({\mathcal{B}},{\mathcal{A}})$).
\end{corollary}

\medskip

We would like to conclude this paper with a discussion of the behaviour of
the set of semi-normalizers (normalizers) with respect to addition. It is
clear that, if $T$ and $S$ are semi-normalizers (normalizers) of an algebra
into another, the sum $T+S$ is not necessarily a semi-normalizer
(normalizer). Proposition \ref{p62} below gives a necessary and sufficient
condition for this to happen, when the algebras in question are CSL
algebras. First we prove a statement, concerning arbitrary masa-bimodules.

\begin{lemma}
\label{lemma60} Let $\Gamma :{\mathcal{B}}({{\mathcal{H}}_{1}},{{\mathcal{H}}%
_{2}})\longrightarrow \{0,1\}$ be such that, if $\Gamma (T)=1$ for some
operator $T$, then there is a reflexive masa-bimodule ${\mathcal{U}}$,
containing $T$, such that $\Gamma (S)=1$ for all $S\in {\mathcal{U}}$.
Suppose that $T_{1}$ and $T_{2}$ are operators, such that $\Gamma (\lambda
_{1}T_{1}+\lambda _{2}T_{2})=1$ for all real numbers $\lambda _{1}$ and $%
\lambda _{2}$. Then there exists a reflexive masa-bimodule ${\mathcal{U}}$,
containing $T_{1}$ and $T_{2}$, such that $\Gamma (T)=1$ for all $T\in {%
\mathcal{U}}$.
\end{lemma}

\proof
Represent ${{\mathcal{H}}_{1}}$ as $L^{2}(X,\mu)$ and ${{\mathcal{H}}_{2}}$
as $L^{2}(Y,\nu)$, where $(X,\mu)$ and $(Y,\nu)$ are compact metric spaces
with regular Borel measures, and let ${\mathcal{D}_{1}}$ and ${\mathcal{D}%
_{2}}$ be the respective multiplication algebras. By the support $%
\mathop{\rm supp}T$ of an operator $T$ we will mean the (closed) support of
the reflexive ${\mathcal{D}_{2}},{\mathcal{D}_{1}}$-bimodule generated by $T$%
. Let $\mathop{\rm supp}T_{1}=\kappa_{1}$, $\mathop{\rm supp}%
T_{2}=\kappa_{2} $ and $\kappa=\kappa_{1}\cup\kappa_{2}$. The sets $%
\kappa_{1},\kappa_{2}$ and $\kappa$ are compact subsets of $X\times Y$. It
suffices to find a real number $\lambda$ such that the operator $%
T_{1}+\lambda T_{2}$ has support $\kappa$. Indeed, in this case, there
exists a reflexive ${\mathcal{D}_{2}},{\mathcal{D}_{1}}$-bimodule ${\mathcal{%
U}}$ containing $T_{1}+\lambda T_{2}$ such that $\Gamma(T)=1$ for all $T\in{%
\mathcal{U}}$. Since $\kappa \subseteq\mathop{\rm supp}\mathcal{U}$ it
follows from \cite{eks}, Theorem 4.6 that $T_{1},T_{2}\in{\mathcal{U}}$.

Put $\kappa _{\lambda }=\mathop{\rm supp}(T_{1}+\lambda T_{2})$, $\upsilon
_{\lambda }=\kappa _{\lambda }^{c}$. It is clear that $\upsilon _{\lambda }$
is an open subset of $X\times Y$. We claim that, if $\lambda \neq \lambda
^{\prime }$ are nonzero, then $\upsilon _{\lambda }\cap \upsilon _{\lambda
^{\prime }}\subseteq \kappa ^{c}$. Indeed whenever $\alpha $ and $\beta $
are open and such that $\alpha \times \beta \subseteq \upsilon _{\lambda
}\cap \upsilon _{\lambda ^{\prime }}$, then $F(\beta )(T_{1}+\lambda
T_{2})E(\alpha )=0$ and also $F(\beta )(T_{1}+\lambda ^{\prime
}T_{2})E(\alpha )=0$. Thus $(\lambda -\lambda ^{\prime })F(\beta
)T_{2}E(\alpha )=0$ and hence $F(\beta )T_{2}E(\alpha )=0$, so that $\alpha
\times \beta \cap \kappa _{2}=\emptyset $. Repeating the argument to the
operators $\frac{1}{\lambda }T_{1}+T_{2}$ and $\frac{1}{\lambda ^{\prime }}%
T_{1}+T_{2}$ gives $\alpha \times \beta \cap \kappa _{1}=\emptyset $. It
follows that $\upsilon _{\lambda }\cap \upsilon _{\lambda ^{\prime }}\cap
\kappa _{2}=\emptyset $ and $\upsilon _{\lambda }\cap \upsilon _{\lambda
^{\prime }}\cap \kappa _{1}=\emptyset $, so that $\upsilon _{\lambda }\cap
\upsilon _{\lambda ^{\prime }}\subseteq \kappa ^{c}$.

Thus the union $\cup\{\upsilon_{\lambda}\cap\kappa:\lambda\in\mathbb{R}%
\backslash\{0\}\}$ is an uncountable union of (relatively) open disjoint
subsets of $\kappa$. Since $\kappa$ is a compact metric space, no more than
countably many of them can be nonempty. Thus there exists $\lambda \in%
\mathbb{R}\backslash\{0\}$ so that $\upsilon_{\lambda}\cap\kappa=\emptyset$
or $\kappa\subseteq\kappa_{\lambda}$. But $\kappa_{\lambda}\subseteq\kappa$.
Indeed if $(\alpha\times\beta)\cap\kappa=\emptyset$ with $\alpha$ and $\beta$
open, then $F(\beta)T_{i}E(\alpha)=0$ for $i=1,2$ so $F(\beta)(T_{1}+\lambda
T_{2})E(\alpha)=0$ which means that $(\alpha\times\beta)\cap\kappa_{\lambda
}=\emptyset$. $\diamondsuit$\hfill\bigskip

\begin{proposition}
\label{p62} Let ${\mathcal{A}}$ and ${\mathcal{B}}$ be CSL algebras and $%
T,S\in SN({\mathcal{B}},{\mathcal{A}})$ (resp. $T,S\in N({\mathcal{B}},{%
\mathcal{A}})$). Then $T+S\in SN({\mathcal{B}},{\mathcal{A}})$ ($T+S\in N({%
\mathcal{B}},{\mathcal{A}})$) if and only if there is a reflexive
normalizing masa-bimodule ${\mathcal{U}}\subseteq SN({\mathcal{B}},{\mathcal{%
A}})$ (${\mathcal{U}}\subseteq N({\mathcal{B}},{\mathcal{A}})$) such that $%
T,S\in {\mathcal{U}}$.
\end{proposition}

\proof
We will consider the case of semi-normalizers, the case of normalizers is
similar. It is clear that if there exists a reflexive normalizing
masa-bimodule ${\mathcal{U}}\subseteq SN({\mathcal{B}},{\mathcal{A}})$ such
that $T,S\in{\mathcal{U}}$, then $T+S$ is a semi-normalizer. Conversely,
suppose that $T,S$ and $T+S$ are semi-normalizers of ${\mathcal{B}}$ into ${%
\mathcal{A}}$. An elementary computation shows that then $%
T^{\ast}BS+S^{\ast}BT\in{\mathcal{A}}$ for each $B\in{\mathcal{B}}$. But
then it is easy to see that $\lambda_{1}T+\lambda_{2}S$ is a semi-normalizer
for all real numbers $\lambda_{1}$ and $\lambda_{2}$. From Theorem \ref{th17}
we have that, if $T\in SN({\mathcal{B}},{\mathcal{A}})$, then there exists a
reflexive space ${\mathcal{U}}\subseteq SN({\mathcal{B}},{\mathcal{A}})$,
containing $T$. The space $\u$ is a masa bimodule, since the semi-lattices
of its map are commutative. Thus Lemma \ref{lemma60} applies with the property of
being a semi-normalizer of ${\mathcal{B}}$ into ${\mathcal{A}}$ in the place
of $\Gamma$. $\diamondsuit$\hfill\bigskip

\begin{corollary}
\label{c64} Let ${\mathcal{A}}$ and ${\mathcal{B}}$ be CSL algebras and $T,S$
be semi-normalizers (normalizers) of ${\mathcal{B}}$ into ${\mathcal{A}}$,
such that $T+S$ is also a semi-normalizer (normalizer). Then $%
B_{1}TA_{1}+B_{2}SA_{2}$ is a semi-normalizer (normalizer) for every $%
B_{1},B_{2}\in {\mathcal{B}}\cap {\mathcal{B}}^{\ast }$ and $A_{1},A_{2}\in {%
\mathcal{A}}\cap {\mathcal{A}}^{\ast }$.
\end{corollary}

\noindent\textbf{Addendum} After this paper was completed, we were
informed that what we call ``normalizing spaces of operators''
have been studied, from a different viewpoint, by other authors
under the names ``ternary rings of operators'' or ``triple
systems''. See for example M. Neal and B. Russo, Contractive
Projections and Operator Spaces (preprint, arXiv: math.OA/0201187)
and its bibliography. We thank J. Arazy, D.P. Blecher and B. Russo
for bringing the relevant literature to our attention. We have
chosen to retain the name ``normalizing spaces" to emphasize the
relation with normalizers of operator algebras, which is one of
the main points of our work.

We have also found out that our Proposition 2.6 is essentially
contained in section 3 of the paper of L.A. Harris, `A
generalization of C*-algebras', Proc. London Math. Soc. (3) 42
(1981) 331-361.

\end{document}